\documentclass[12pt,a4paper,reqno]{article}
\def\hybrid{\topmargin 0pt      \oddsidemargin 0pt
        \headheight 0pt \headsep 0pt
        \textwidth 160true mm       % US paper
        \textheight 231true mm         % US paper
        \marginparwidth 0.0in
        \parskip 0pt plus 1pt   \jot = 1.5ex}
\usepackage{amssymb}
\usepackage{amsmath}
\usepackage{amsthm}
\hybrid

\newcommand{\g}{\mathfrak{g}}
\newcommand{\hh}{\mathfrak{h}}
\newcommand{\mm}{\mathfrak{m}}
\renewcommand{\aa}{\alpha}
\newcommand{\ff}{\varphi}

\newcommand{\tb}{\mathbf{t}}
\newcommand{\Uh}{{U}_h (\g)}
\newcommand{\Ug}{U(\g)}

\newcommand{\A}{\mathcal{A}}

\newcommand{\C}{\mathbb{C}}
\newcommand{\D}{\mathcal{D}}
\renewcommand{\O}{\mathcal{O}}

\newcommand{\Z}{\mathcal{Z}}
\newcommand{\X}{\mathcal{X}}

\newcommand{\I}{\mathcal{I}}
\newcommand{\ot}{\otimes}

\renewcommand{\[}{[\![}
\renewcommand{\]}{]\!]}

\newcommand{\hlf}[1]{\frac{#1}{2}}

\newcommand{\qqquad}{\quad\quad\quad}

\newcommand{\ad}{\operatorname{ad}}
\newcommand{\id}{\operatorname{id}}
\newcommand{\Id}{\operatorname{Id}}
\newcommand{\Tr}{\operatorname{Tr}}
\newcommand{\End}{\operatorname{End}}
\newcommand{\bL}{\overline{L}}

\newcommand{\vep}{\varepsilon}
\renewcommand{\Im}{\operatorname{Im}}
\newcommand{\Ker}{\operatorname{Ker}}

\renewcommand{\ll}{\lambda}

\newcommand{\bsig}{\bar{\sigma}_h}
\newcommand{\Ga}{\Gamma}
\newcommand{\bO}{\overline\Omega_\Gamma}

\theoremstyle{plain} %% This is the default

\newtheorem{cor}{Corollary}[section]
\newtheorem{lemma}{Lemma}[section]
\newtheorem{propn}{Proposition}[section]

\theoremstyle{definition}
\newtheorem{defn}{Definition}[section]

\newtheorem{re}{}[section]
\newtheorem{conjecture}[re]{Question}

\theoremstyle{definition}  %{remark}
\newtheorem{rem}{Remark}[section]

\newcommand{\be}[1]{\begin{eqnarray#1}}
\newcommand{\ee}[1]{\end{eqnarray#1}}

\newcommand{\tor}[1]{\stackrel{#1}{\longrightarrow}}

\newcommand{\De}{\Delta}
\newcommand{\de}{\delta}

\newcommand{\cA}{\mathcal{C}}
\renewcommand{\t}{\otimes}
\newcommand{\wt}{\widetilde}

\newcommand{\bbe}{{\bar\beta}}

\newcommand{\Ea}{E_{\aa}}
\newcommand{\Ema}{E_{-\aa}}
\newcommand{\rc}{\Omega^+\setminus\Omega_\Ga}
\newcommand{\baa}{{\bar\aa}}

\numberwithin{equation}{section}

\begin{document}
\newcommand{\kh}{\C[[h]]}
\newcommand{\dl}{[\![}
\newcommand{\dr}{]\!]}
\newcommand{\Vect}{\operatorname{Vect}}
\newcommand{\Ze}{{\mathbb Z}}
\title{Double quantization on coadjoint representations of simple Lie groups
and its orbits}
                        
\author
{Joseph Donin\\
{\normalsize Max-Planck-Institut f\"ur Mathematik}}

\date{}
\maketitle

\begin{abstract}
Let $M$ be a manifold with an action of a Lie group $G$, $\A$ the function
algebra on $M$. 
The first problem we consider is to construct a $U_h(\g)$ invariant quantization, $\A_h$,
of $\A$, where $U_h(\g)$ is a quantum group corresponding to $G$.

Let $s$ be a $G$ invariant Poisson bracket on $M$. The second problem we
consider is to construct a $U_h(\g)$ invariant two parameter (double) quantization, $\A_{t,h}$,
of $\A$ such that $\A_{t,0}$ is a $G$ invariant quantization of $s$.
We call $\A_{t,h}$ a $U_h(\g)$ invariant quantization of the Poisson bracket $s$.

In the paper we study the cases when $G$ is a simple Lie group and $M$
is the coadjoint representation $\g^*$ of $G$ or a semisimple orbit in this representation.

First of all, we describe Poisson brackets and pairs of Poisson brackets
related to $U_h(\g)$ invariant quantizations for arbitrary algebras. After that
we construct a two parameter quantization on $\g^*$
for $\g=sl(n)$ and $s$ the Lie bracket and show that such a quantization does not exist for other
simple Lie algebras. 
As the function algebra on $\g^*$ we take the symmetric algebra $S\g$.
In $sl(n)$ case, we also consider the problem of restriction of  
the family $(S\g)_{t,h}$ on orbits. In particular, we describe explicitly
the Poisson bracket along the parameter $h$ of this family, which turns out to be
quadratic, and prove that it can be restricted on each orbit in $\g^*$.
We prove also that the family $(S\g)_{t,h}$ can be restricted on the maximal semisimple orbits.

For $M$ a manifold isomorphic to a semisimple orbit in $\g^*$, 
we describe the variety of all brackets
related to the one parameter quantization. Actually, it is a variety making $M$   
into a Poisson manifold with a Poisson action of $G$. 
It turns out that not all such brackets and not all orbits admit a double quantization
with $s$ the Kirillov-Kostant-Souriau bracket. We classify the orbits and pairs of brackets
admitting a double quantization and construct such a quantization for almost all
admissible paires.

\end{abstract}

\section{Introduction}

Quantum groups can be considered as symmetry objects
of certain ``quantum spaces'' described by noncommutative
algebra of functions. This point of view was developed, for example,
in \cite{RTF} and \cite{Ma}.
Here we study the inverse problem: given the quantum group corresponding
to a Lie group $G$,
we want to define a ``quantum space'' corresponding to
a given classical $G$-manifold.  

Let $M$ be a manifold with an action of a  Lie group $G$,
$\g$ the Lie algebra of $G$, and $\Uh$ the quantized universal enveloping
algebra. 
Let $\A$ be the sheaf of function algebras on $M$. It may be a sheaf
of smooth, analytic, or algebraic functions. For shortness,
we simply call $\A$ a function algebra. The algebra $\A$ is of course
invariant under the induced action of the bialgebra $U(\g)$.

We consider the following two general problems.

{\bf The first problem.} Does there exist a deformation quantization,
$\A_h$, of $\A$, which is invariant under the action of the
quantum group $U_h(\g)$?

{\bf The second problem.} Suppose $\A_t$ is a $U(\g)$ invariant 
quantization of
$\A$. Does there exist a two parameter quantization, $\A_{t,h}$, of
$\A$ such that $\A_{t,0}=\A_t$, which is invariant under $\Uh$?  

In this paper, we study the first and the second problems for two cases.
The first case, when $M$ is the 
coadjoint representation of a simple Lie group. The second case, 
when $M$ is a semisimple orbit
in this representation. 
This paper is motivated by papers \cite{Do2} and \cite{DGS} where
we started to study these problems. In this paper
we develop  results of \cite{Do2} and \cite{DGS} and  
present some additional results.

The paper is organized as follows.

In Section 2 we recall some facts about quantum groups and
related categories, 
which are essential for a strict formulation of our problems
and for our approach to 
$\Uh$ invariant quantization of algebras.
In particular, we use the Drinfeld category with non-trivial associativity 
constraint determined by
an invariant element $\Phi_h\in\Ug^{\ot 3}[[h]]$ and show
that the problem of $\Uh$ invariant quantization is equivalent to
the problem of deforming the function algebra in such a way that
the deformed algebra to be $G$ invariant and $\Phi_h$ associative
(see Subsection 2.3).

Subsection 2.4 is very important for the paper. In this subsection we give, for all 
commutative algebras,
a description of Poisson brackets related to one and two parameter
$\Uh$ invariant quantizations.
We show the following.
If $\A_h$ is a
$U_h(\g)$ invariant quantization, the corresponding Poisson bracket, $p$,
on $M$ has to be a difference of two brackets, $p=f-r_M $. Here
$r_M$ is the so called $r$-matrix bracket obtained from a classical
$r$-matrix $r\in \wedge^2\g$ with the help of the action morphism
$\g\to \Vect(M)$. So, the Schouten bracket $\dl r_M,r_M \dr$ is equal to
the image $\ff_M$ of the invariant element $\ff\in \wedge^3\g$.
The bracket $f$ is $U(\g)$ invariant and such that $\dl f,f \dr=-\ff_M$. 
Of course, any invariant bracket, $f$, is compatible with $r_M$, so that
$\dl p,p \dr=0$. 

We see that for existence of the family $\A_h$ one needs  existence
of an invariant bracket $f$ on $M$ such that 
$$\dl f,f \dr=-\ff_M. \eqno(1.1)$$
Note that the manifold $M$ endowed with the bracket $p=f-r_M$
is a Poisson manifold with a Poisson action of $G$, where $G$ is considered
to be the Poisson-Lie group with Poisson structure defined by $r$.
We shall not use this fact in the paper. 

Similarly, given a two parameter quantization, $\A_{t,h}$, a pair
of compatible Poisson brackets is determined. 
These brackets are: the bracket $p=f-r_M$ considered above
and a $U(\g)$ invariant Poisson bracket, $s$, the initial term
of the $U(\g)$ invariant quantization $\A_t$. 
We may perceive the family $\A_{t,h}$ as a $\Uh$ invariant
quantization of the Poisson bracket $s$. 

We assume that $s$ is given
in advance 
and determined, for example, by a $G$ invariant simplectic structure on $M$.
From the compatibility of $p$ and $s$ (this means $\dl p,s \dr=0$) 
follows that 
$$\dl f,s \dr=0. \eqno(1.2)$$

So, for existence of the family $\A_{t,h}$ one needs existence
of an invariant bracket $f$ on $M$ such the both equations (1.1)
and (1.2) hold. 

Thus, our problems divide into two steps. The first step is looking for
invariant brackets $f$ on $M$ satisfying either (1.1) (in case of the first problem)
or both (1.1) and (1.2) (in case of the second problem).
The second step is quantizing these brackets.

In Section 3 we consider the one and two parameter quantization on $M=\g^*$,
the coadjoint representation of a simple Lie algebra $\g$.
As a function algebra on $\g^*$, we take the symmetric algebra $S\g$.
It turns out that the cases $\g=sl(n)$ and $\g\neq sl(n)$ are
quite different. 

We prove that for $\g\neq sl(n)$ the two parameter
family which is a $\Uh$ invariant quantization of the Lie
bracket on $S\g$ does not exist. Moreover, as a conjecture we state
that in this case even a one parameter $\Uh$ invariant quantization of $S\g$
does not exist. 

In the case $\g=sl(n)$, the two parameter quantization of $S\g$  exists.
Moreover, the picture looks like in the classical case.
Recall that in the classical case, the natural one parameter
$\Ug$ invariant quantization of $S\g$ is given by the family
$(S\g)_t=T(\g)[t]/J_t$, where $J_t$ is  the ideal generated by the elements
of the form $x\ot y-\sigma(x\ot y)-t[x,y]$, $x,y\in\g$, $\sigma$ is
the permutation. By the PBW theorem, $(S\g)_t$ is a free module
over $\C[t]$. 
We have $(S\g)_0=S\g$, so this
family of quadratic-linear algebras gives a $U(\g)$ invariant 
quantization of $S\g$. It is obvious that the Poisson
bracket, $s$, related to this quantization is the Lie bracket on $\g^*$.

We show that for $\g=sl(n)$ this picture can be extended to 
the quantum case.
Namely, there exist deformations, $\sigma_h$ and $[\cdot,\cdot]_h$, of 
both the mappings
$\sigma$ and $[\cdot,\cdot]$ such that
the two parameter family of algebras
$(S\g)_{t,h}=T(\g)[[h]][t]/J_{t,h}$, where $J_{t,h}$ is the ideal generated 
by the elements
of the form $x\ot y-\sigma_h(x\ot y)-t[x,y]_h$, $x,y\in\g$,
gives a $\Uh$ invariant quantization of the Lie bracket $s$ on $\g^*$.
In this case, the corresponding bracket $f$ from (1.2) is a quadratic
bracket which is, up to a factor, a unique nontrivial invariant map 
$\wedge^2\g\to S^2\g$.

Taking $t=0$ we obtain the family $(S\g)_h$ which is a quadratic algebra
over $\C[[h]]$. This algebra can be called the quantum symmetric algebra 
(or quantum polynomial
algebra on $\g^*$). We show (Subsection 3.4) that $(S\g)_h$ can be
included in the deformed graded differential algebra (deformed de Rham complex).
In Subsection 3.5 we prove that the family $(S\g)_{t,h}$ can be 
restricted on the maximal semisimple orbits in $\g^*$ to
give a two parameter quantization on these orbits.

In Section 4 we study the problems of one and two parameter quantization
on semisimple orbits in $\g^*$ for all simple Lie algebras $\g$.
First of all, we classify all the brackets $f$ satisfying (1.1) and both (1.1) and (1.2) for $s$
being the Kirillov-$\A_{t,h}$ (KKS) bracket on the orbit.
After that, we construct quantizations of these brackets.

Let $M$ be a semisimple orbit. In Subsection 4.1 we prove that the brackets $f$ satisfying
(1.1) form a $\dim H^2(M)$-dimensional variety. We give a description
of this variety and prove (in Subsection 4.3)  that almost
all these brackets can be quantized. So, we obtain for $M$ 
a $\dim H^2(M)$ parameter family of non-equivalent one parameter quantizations.
Note that in \cite{DG2} we have built one of these quantizations, the quantization
of the so called Sklyanin-Drinfeld Poisson bracket.

It turns out that brackets $f$ satisfying (1.1) and (1.2) exist not
for all orbits.
We call an orbit $M$ {\em good} if there exists a bracket $f$ satisfying 
(1.1) and (1.2) 
for the Kirillov-Kostant-Souriau (KKS) bracket $s$.

In Subsection 4.1 we give  
the following classification of the semisimple good orbits for all 
simple $\g$, \cite{DGS}.

In the case $\g=sl(n)$ all semisimple orbits are good.
(Actually we prove that in this case all orbits are good.)

For $\g\neq sl(n)$ all symmetric orbits (which are symmetric spaces) 
are good. In this case $\ff_M=0$,
so $r_M$ itself is a Poisson bracket compatible with $s$. 

Only in the case $\g$ of type $D_n$ and $E_6$ (except of $A_n$) there
are good orbits different from the symmetric ones. For such orbits
$\ff_M\neq 0$. 

We show that brackets $f$  on a good orbit satisfying (1.1) and (1.2),
form a one parameter family.

In Subsection 4.2 we consider cohomologies of an invariant complex with the differential
given by the Schouten bracket with the bivector $f$.  These cohomologies are
needed for our construction of quantization.

In Subsection 4.3 we construct one and two parameter quantizations
for semisimple orbits. 
According to our approach, as a first step 
we construct a $G$ invariant  
$\Phi_h$ associative quantization, i.e., a quantization in the Drinfeld category
with  non-trivial associativity constraint given by $\Phi_h$.
Note that the bracket $f$ from (1.1) can be considered as a ``Poisson bracket''
in that category.
As a second step, we make a passage to the category with trivial associativity
to obtain the associative $\Uh$ invariant quantization.
We applied this method earlier for quantizing the function algebra on
the highest weight orbits in irreducible representations of $G$, the algebra of sections
of linear vector bundles over flag manifolds, and the function algebra
on symmetric spaces, \cite{DGM}, \cite{DG1}, \cite{DS1}.  

I put in the text some questions
which  naturally appeared by exposition.
They are open (for me) and seem to be important.

This research is partially supported
by Israeli Academy Grant N 8007/99 
and by the Emmy Noether Institute and the Minerva
Foundation of Germany

{\bf Acknowledgment.}
I am very grateful to J.Bernstein, D.Gurevich, S.Shnider, and V.Ostapenko
for many useful discussions.

I thank Max-Planck-Institut f\"ur Mathematik for hospitality
and very stimulating working atmosphere.

\section{Preliminaries}
\subsection{Quantum groups}
We shall consider quantum groups in sense of Drinfeld, \cite{Dr2},
as deformed universal enveloping algebras.  
If $U(\g)$ is the universal enveloping algebra of a complex Lie algebra  $\g$, then
the quantum group (or quantized universal enveloping algebra) corresponding
to $U(\g)$ is a topological Hopf algebra, $\Uh$, over $\kh$, isomorphic to $U(\g)[[h]]$ as a
topological $\kh$ module and such that $\Uh/h\Uh=U(\g)$ as a Hopf algebra
over $\C$.
In particular, the deformed comultiplication in  $\Uh$ has the form
\be{}
\label{comult}
\De_h=\De+h\De_1+o(h),
\ee{}
where $\De$ is the comultiplication in the universal enveloping algebra $U(\g)$.
One can prove, \cite{Dr2},  that
$\De_1:U(\g)\to U(\g)\ot U(\g)$ is such a map that $\De_1-\sigma\De_1=\de$
($\sigma$ is the usual permutation) being restricted on $\g$ gives
a map $\de:\g\to\wedge^2\g$ which is a 1-cocycle and defines
the structure of a Lie coalgebra on $\g$ (the structure of a Lie algebra on
the dual space $\g^*$).  
The pair $(\g,\de)$ is
considered as a quasiclassical limit of $\Uh$.

In general, a pair $(\g,\de)$, where
$\g$ is a Lie algebra and $\de$ is such a 1-cocycle, is called a Lie bialgebra.
It is proven, \cite{EK}, that any Lie bialgebra $(\g,\de)$ can be quantized, i.e., 
there exists a quantum group  $\Uh$ such
that  the  pair  $(\g,\de)$ is its quasiclassical limit.

A Lie bialgebra $(\g,\de)$ is said to be  a coboundary one if  there exists an element
$r\in\wedge^2$, called the classical $r$-matrix, 
such that $\de(x)=[r,\De(x)]$ for $x\in\g$. Since $\de$ defines a Lie coalgebra
structure, 
$r$ has to satisfy the so-called classical Yang-Baxter equation
which can by written in the form
\be{}
\label{cYB}
\[r,r\]=\ff,
\ee{}
where $\[\cdot,\cdot\]$ stands for the Schouten bracket and 
$\ff\in\wedge^3\g$ is an invariant element. We denote the coboundary
Lie bialgebra by $(\g,r)$.

In case $\g$ is a simple Lie algebra, the most known is
the Sklyanin-Drinfeld $r$-matrix:
\be{*}
r=\sum_\alpha X_\alpha\wedge X_{-\alpha},
\ee{*} 
where the sum runs over all positive roots; the root vectors $X_\alpha$
are chosen is such a way that $(X_\alpha, X_{-\alpha})=1$
for the Killing form $(\cdot,\cdot)$.
This is the only $r$-matrix of weight zero, \cite{SS}, and its quantization is the Drinfeld-Jimbo
quantum group.
A classification of all $r$-matrices for simple Lie algebras was
given in \cite{BD}.

We are interested in the case when $\g$ is a semisimple finite dimensional Lie algebra.
In this case, from results of Drinfeld and Etingof and Kazhdan one can derive
the following
\begin{propn}
\label{prop2.1}
Let $\g$ be a semisimple Lie algebra. Then 

a) any Lie bialgebra $(\g,\de)$ is a coboundary one;

b) the quantization, $\Uh$, of any coboundary Lie bialgebra $(\g,r)$ exists and
is isomorphic to $\Ug[[h]]$ as a topological $\C[[h]]$ algebra;

c) the comultiplication in $\Uh$ has the form
\be{}
\label{comul}
\De_h(x)=F_h\De(x)F^{-1}_h, \qquad x\in\Ug,
\ee{}
 where $F_h\in \Ug^{\ot 2}[[h]]$ and can be chosen in the form
\be{}
\label{F}
F_h=1\ot 1+\hlf{h}r +o(h).
\ee{} 
\end{propn}

\begin{proof}
a) follows from the fact that $H^1(\g,\wedge^2\g)=0$.
From the fact that $H^2(\g,\Ug)=0$ follows that
$\Ug$ does not admit any
nontrivial deformations as an algebra, (see \cite{Dr1}), which
proves b). 
From the fact that $H^1(\g,\Ug^{\ot 2})=0$ follows that
any deformation of the algebra morphism
$\De:\Ug\to\Ug\ot\Ug$ appears as a conjugation of $\De$.
In particular, the comultiplication in $\Uh$ looks like (\ref{comul}) 
with some $F_h$ such that $F_0=1\ot 1$.

From the coassociativity of $\De_h$ follows that $F_h$  
satisfies  the
equation
\be{}
\label{eqF}
(F_h\t 1)\cdot(\De\t id)(F_h)=(1\t F_h)\cdot(id\t\De)(F_h)\cdot\Phi_h
\ee{}
for some invariant element $\Phi_h\in\Ug^{\ot 3}[[h]]$.

The element $F_h$ satisfying (\ref{comul}) and (\ref{F}) can be obtained by correction
of some $F_h$ only obeying (\ref{comul}), \cite{Dr2}.
This procedure also makes use simple cohomological arguments 
and essentially (\ref{eqF}). This proves c).
\end{proof}

From (\ref{eqF}) follows that if $F_h$ has the form (\ref{F}), then
the coefficient by $h$ for $\Phi_h$ vanishes. Moreover, as a coefficient by $h^2$
one can take the element $\ff$ from (\ref{cYB}), i.e.,  
\be{}\label{fPhi}
\Phi_h=1\t 1\t 1+h^2\ff+o(h^2).
\ee{}
In addition, from (\ref{eqF}) follows that
$\Phi_h$ satisfies the pentagon identity
\be{}
\label{pent}
(id^{\t 2}\t \De)(\Phi_h)\cdot(\De\t id^{\t 2})(\Phi_h)=
(1\t\Phi_h)\cdot(id\t\De\t id)(\Phi_h)\cdot(\Phi_h\t 1).
\ee{}

\begin{conjecture}\label{con2.1}
Let $(\g,r)$ be a coboundary Lie bialgebra.
Does there exist a quantization of it, $\Uh$, such that
$\Uh$ is isomorphic to $\Ug[[h]]$ as a topological $\C[[h]]$ algebra  
and the comultiplication has the form (\ref{comul})?
\end{conjecture}

From \cite{Dr4} follows that if $\[r,r\]=0$, the answer to this
question is positive.

\subsection{Categorical interpretation}

It is known that the elements constructed above 
have a nice
categorical interpretation.
First, recall some facts about the Drinfeld algebras and the monoidal
categories determined by them.

Let $A$ be a commutative algebra with unit, $B$ a unitary
$A$-algebra. The category of representations of $B$ in $A$-modules,
i.e. the category of $B$-modules, will be a monoidal category 
if the algebra $B$ is equipped with an algebra morphism, $\De: B\to B\ot_A B$, 
 called  comultiplication,  and an invertible element
$\Phi\in B^{\t 3}$ such that $\De$ and $\Phi$ 
satisfy the conditions (see \cite{Dr2})
\be{}
&(id\t\De)(\De(b))\cdot\Phi=\Phi\cdot (\De\t id)(\De(b)),\ \ b\in B, 
\label{d1}\\
&(id^{\t 2}\t \De)(\Phi)\cdot(\De\t id^{\t 2})(\Phi)=
(1\t\Phi)\cdot(id\t\De\t id)(\Phi)\cdot(\Phi\t 1).  \label{d2}       
\ee{}
Define a tensor product functor
for  the category of
$B$ modules $\cA$,  denoted  $\ot_{\cA}$  or simply $\t$ when 
there can be no confusion, in
the following way: given $B$-modules $M,N$,
$M\t_{\cA} N=M\t_A N$ as an $A$-module. The
action of $B$ is defined by
 $$b(m\t n)=(\De b)(m\t n)= b_1m\t b_2n,$$
where $\De b=b_1\t b_2$
(we use the Sweedler convention of an implicit summation over an index).
The element $\Phi=\Phi_1\t\Phi_2\t\Phi_3$ defines the
 associativity constraint,
$$a_{M,N,P}:(M\t N)\t P\to M\t(N\t P),\,\, 
a_{M,N,P}((m\t n)\t p)= \Phi_1m\t(\Phi_2n\t\Phi_3p).$$
Again the summation in the expression for $\Phi$ is understood.
By virtue of (\ref{d1}) $\Phi$ induces an isomorphism of $B$-modules, and by
virtue of
(\ref{d2}) the pentagon identity for monoidal categories holds. We call the
triple
$(B,\De,\Phi)$ a Drinfeld algebra.
The definition is somewhat non-standard in that
 we do not require the existence of an antipode.
The category $\cA$ of $B$-modules for $B$ a 
Drinfeld algebra  becomes a monoidal category. When
it becomes necessary to be more explicit  we shall denote 
${\cA}(B,\De,\Phi)$. 

Let $(B,\De,\Phi)$ be a Drinfeld algebra and $F\in B^{\t 2}$
an invertible element. Put 
\be{}\label{f1}
&\wt{\De}(b)=F\De(b)F^{-1},\ \ b\in B,  \\
&\wt{\Phi}=(1\t F)\cdot(id\t\De)(F)\cdot\Phi\cdot(\De\t id)(F^{-1})
\cdot(F\t 1)^{-1}.           \label{f2}
\ee{}
Then $\wt{\De}$ and $\wt{\Phi}$ satisfy (\ref{d1}) and (\ref{d2}),
therefore the triple $(B,\wt{\De},\wt{\Phi})$ also becomes a Drinfeld
algebra. We say that it is obtained by twisting from   $(B,\De,\Phi)$.
It has an equivalent monoidal category of modules, 
$\wt{\cA}(B,\wt{\De},\wt{\Phi})$. Note that the equivalent categories $\cA$
and $\wt{\cA}$ consist of the same objects as $B$-modules,
and the tensor products of two objects are isomorphic as $A$-modules.
 The equivalence $\cA\to \wt{\cA}$ is given by the pair $(\Id,F)$,
where $\Id:\cA\to \wt{\cA}$ is the identity functor of the categories
(considered without the monoidal structures, but only as categories
of $B$-modules), and $F:M\t_{\cA} N\to M\t_{\wt{\cA}} N$ is defined by 
$m\t n\mapsto F_1m\t F_2n$ where $F_1\t F_2=F$.

We are interested in the case when $A=\C[[h]]$, $B=U(\g)[[h]]$
where $\g$ is a complex semisimple Lie algebra. In this case, all
tensor products over $\C[[h]]$ are completed in
$h$-adic topology.

We have two nontrivial Drinfeld algebras. The first is $(U(\g)[[h]],\De,\Phi_h)$,
with the usual comultiplication and $\Phi_h$ from (\ref{eqF}).
The condition (\ref{d1}) means the invariantness of $\Phi_h$, while (\ref{d2})
coincides with (\ref{pent}). 
The second Drinfeld algebra is $(U(\g)[[h]],\De_h, {\bf 1})$.  
It obtaines by twisting of the first one
by the element $F_h$ from (\ref{comul}). The equation (\ref{f2})
follows from (\ref{eqF}). 
The pair $(\Id,F_h)$ defines an
equivalence between the  corresponding monoidal categories 
${\cA}(U(\g)[[h]],\De,\Phi_h)$
and ${\cA}(U(\g)[[h]],\De_h,{\bf 1})$. The last is the category of
representations of the quantum group $\Uh$.

It is clear that reduction modulo $h$
defines a functor from  either of these categories
to the category of representations of $U(\g)$ and the equivalence
just described reduces to the identity modulo $h$. 
In fact, both categories are 
$\C[[h]]$-linear extensions (or deformations) of the $\C$-linear
category of representations of $\g$. Ignoring the monoidal structure
the extension is a trivial one, but the associator $\Phi_h$ in the first case
and the  comultiplication $\De_h$ in the second case make the
extension non-trivial from the point of view of monoidal categories.

\subsection{$\Uh$ ivariant quantizations of algebras}
\label{ss2.3}
Let  $(B,\De,\Phi)$ be a Drinfeld algebra. 
Assume $\A$ is a $B$-module with a multiplication $\mu:\A\t_A \A\to \A$
which is a homomorphism of $A$-modules. We say that $\mu$
is $\De$ invariant if 
\be{}\label{finv}
b\mu(x\otimes y)=\mu\De(b)(x\t y)\ \ \ \mbox{for}\ b\in B,\ x,y\in \A,
\ee{}
and $\mu$ is $\Phi$  associative, if
\be{}\label{fass}
\mu(\Phi_1x\otimes \mu(\Phi_2 y\t \Phi_3 z)))=\mu(\mu(x\t y)\t z)
\ \ \ \mbox{ for } x,y,z\in \A. 
\ee{}

Note, that a $B$-module $\A$ equipped with $\De$ invariant and
$\Phi$ associative multiplication is an associative algebra in
the monoidal category ${\cA}(B,\De,\Phi)$.
If $(B,\wt{\De},\wt{\Phi})$ is a Drinfeld algebra twisted 
by (\ref{f1}) and (\ref{f2}), then the algebra $\A$
may be transfered into the equivalent category $\wt{\cA}(B,\wt{\De},\wt{\Phi})$:
the multiplication $\wt{\mu}=\mu F^{-1}:M\t_A M\to M$ is
$\wt{\Phi}$-associative and invariant in the category $\wt{\cA}$.

Let $\A$ be a $U(\g)$ invariant associative algebra, i.e., an algebra with 
$U(\g)$ invariant multiplication $\mu$ in sense
of (\ref{finv}). A deformation (or quantization) of $\A$ is
an associative algebra,  $\A_h$, which is 
isomorphic to $\A[[h]]=\A\t\C[[h]]$ (completed tensor
product) as a $\C[[h]]$-module, with multiplication in
$\A_h$ having  the form
$\mu_h=\mu+h\mu_1+o(h)$.
The algebra $\Ug[[h]]$ is clearly acts on the $\C[[h]]$ module $\A_h$.

We will study  quantizations of $\A$ which will be invariant under
the comultiplication $\De_h$. In other words, $\A_h$ will be an algebra
in the category of representations of the quantum group $\Uh$.
It is clear from the previous Subsection that if $\A_h$ is such a quantization, then
the multiplication $\mu_hF_h$ makes the module $\A[[h]]$
into an algebra in the category ${\cA}(U(\g)[[h]],\De,\Phi_h)$, i.e.,
this multiplication is $\Ug$ invariant and $\Phi_h$ associative. 

We shall see that often it is easier to constract 
$\Ug$ invariant and $\Phi_h$ associative quantization of $\A$.
After that, the ivariant quantization with respect to any quantum
group from Proposition \ref{prop2.1} can be obtained by twisting by 
the appropriate $F_h$.

As an algebra $\A$ we may  take an algebra $\A_t$ that is itself 
a $\Ug$ invariant quantization of a commutative algebra $\A$.
In this case, a $\Uh$ invariant quantization of $\A_t$ is an
algebra $\A_{t,h}$ over $\C[[t,h]]$.

\subsection{Poisson brackets associated with the $\Uh$ invariant\\ quantization}

Let $\A$ be a $U(\g)$ invariant commutative algebra with 
 multiplication $\mu$ and
 $\A_h$ its quantization with multiplication
$\mu_h=\mu+h\mu_1+o(h)$.
The Poisson bracket corresponding to the quantization is given by
$\{a,b\}=\mu_1(a,b)-\mu_1(b,a)$, $a,b\in\A$.

In general, we call a skew-symmetric bilinear form $\A\t\A\to\A$ a bracket, 
if it satisfies the Leibniz rule in either argument when the other is fixed.
The term Poisson bracket indicates that the Jacobi identity is also true.

A bracket of the form 
\be{}\label{rmb}
\{a,b\}_r=(r_1a)(r_2b)=\mu r(a\ot b) \qqquad a,b\in\A, 
\ee{}
where $r=r_1\t r_2$ (summation implicit) is the representation of 
$r$-matrix $r$, will be called an $r$-matrix bracket.

Assume $\A_h$ is a $\Uh$ invariant quantization, i.e., the multiplicatin
$\mu_h$ is $\De_h$ invariant. 
We shall show that in this case the Poisson bracket $\{\cdot,\cdot\}$ has
a special form.  Suppose $f$ and $g$ are two brackets
on $\A$. Define their Schouten bracket $\[f,g\]$ as 
\be{}
\[f,g\](a,b,c)=f(g(a,b),c)+g(f(a,b),c)+{\rm cyclic\ permutations\ of}\ a,b,c.
\ee{}
Then $\[f,g\]$ is a  skew-symmetric map $\A^{\t 3}\to \A$.
We call $f$ and $g$ compatible if $\[f,g\]=0$.

\begin{propn}\label{prop2.2}
Let $\A$ be a $U(\g)$ invariant commutative algebra and
$\A_h$ a $\Uh$ invariant quantization.
Then the corresponding Poisson bracket has the form
\be{}\label{rinv}
\{a,b\}=f(a,b)-\{a,b\}_r
\ee{}
where $f(a,b)$ is a $U(\g)$ invariant bracket.

The brackets $f$ and $\{\cdot,\cdot\}_r$ are
compatible and $\[f,f\]=-\ff_\A$,
where $\ff_\A(a,b,c)=(\ff_1a)(\ff_2b)(\ff_3c)$ 
and $\ff_1\t\ff_2\t\ff_3=\ff\in\wedge^3\g$ is the invariant element
from (\ref{cYB}).
\end{propn}

\begin{proof}
Let the comultiplication for  $\Uh$ have the form (\ref{comult}).
Let $\A$ be a commutative algebra with the $U(\g)$ invariant multiplication $\mu$.
Suppose 
$\A_h$ is a $\Uh$ invariant quantization of $\A$ . 
This means that the deformed multiplication has the form
\be{}
\label{multh}
\mu_h=\mu+h\mu_1+o(h)
\ee{}
and satisfies the relation
\be{}
\label{invh}
x\mu_h(a\ot b)=\mu_h\De_h(x)(a\ot b) \qqquad {\rm for} \qquad x\in U(\g), \ a,b\in \A.
\ee{}

Substituting (\ref{comult}) and (\ref{multh}) in (\ref{invh})
and collecting the terms by $h$ we obtain
$$\mu_1(a\t b)=\mu\De(x)(a\t b)+m\De_1(x)(a\t b).$$
Subtracting from this equation the similar one with permuting $a$ and $b$ 
and making use that $\De$ is commutative
and $\de=\De_1-\sigma\De_1$ is skew-commutative, we derive that the Poisson
bracket $p=\{\cdot,\cdot\}$ has to satisfy the property
\be{}
\label{qpb}
 x p(a\ot b)=p \De(x)(a\ot b) + \mu\de(x)(a\ot b),   \qquad x\in \Ug.
\ee{}  

 Let us prove that the bracket $f(a,b)=\{a,b\}+\{a,b\}_r$ is  $U(\g)$ invariant.
Indeed,  from (\ref{rmb}) 
we have for $x\in U(\g)$, $a,b\in \A$
\be{*}
x\mu r(a\ot b)=\mu\De(x)r(a\ot b)=\mu r\De(x)(a\ot b) - \mu [r,\De(x)](a\ot b).
\ee{*}
Using this expression, (\ref{qpb}), and the fact that $\de(x)=[r,\De(x)]$, we obtain
\be{*}
xf=xp+x\mu r=(p\De(x)+\mu [r,\De(x)]) + (\mu r\De(x)-\mu [r,\De(x)])=\\
=p\De(x)+\mu r\De(x)=f\De(x),
\ee{*}
which proves the invariantness of $f$.

So, we have
$\{a,b\}=f(a,b)-\{a,b\}_r$, as required.

It is easy to check that any bracket of the form
$\{a,b\}=(X_1a)(X_2b)=\mu(X_1a,X_2b)$, for $X_1\t X_2\in \g\wedge \g$, 
 is compatible with any invariant bracket. 
In particular, an $r$-matrix bracket is compatible with $f$.
In addition, $\{\cdot,\cdot\}$ is a Poisson bracket, so its Schouten
bracket with itself is equal to zero. Using this and 
the fact that the Schouten bracket of $r$-matrix bracket with
itself is equal to $\ff_\A$, we obtain from (\ref{rinv}) that
$\[f,f\]=-\ff_\A$.
\end{proof}

\begin{rem}\label{remPL}
Let $\A$ be the function algebra on a $G$-manifold $M$,
where the Lie group $G$ corresponds to the Lie algebra $\g$.
It is easy to see that condition (\ref{qpb}) with $\de(x)=[r,\De(x)]$ is equivalent
to the condition that the pair $(M,p)$ becomes a $(G,\tilde{r})$-Poisson manifold,
where $\tilde{r}$ is the Poisson structure on $G$ defined by the $r$-matrix $r$:
$\tilde{r}=r^\prime-r^{\prime\prime}$, where $r^\prime$
and $r^{\prime\prime}$ are the left- and right-invariant bivector fields on $G$
corresponding to $r$. It is known that $\tilde{r}$ makes $G$ into a Poisson-Lie group.
So Proposition \ref{prop2.2} gives a description of Poisson
structures $p$ on $M$ making $(M,p)$ into a $(G,\tilde{r})$-Poisson manifold.  
\end{rem}

We shall also consider two parameter quantizations of algebras.
A two parameter quantization  of an algebra $\A$ is
an algebra $\A_{t,h}$ isomorphic to $\A[[t,h]]$
as a $\C[[t,h]]$ module and having a multiplication in the form
$$\mu_{t,h}=\mu+t\mu^\prime_1+h\mu^{\prime\prime}_1+o(t,h).$$
With such a quantization, one associates two Poisson brackets:
the bracket $s(a,b)=\mu^\prime_1(a,b)-\mu^\prime_1(b,a)$ along $t$,
and the bracket
$p(a,b)=\mu^{\prime\prime}_1(a,b)-\mu^{\prime\prime}_1(b,a)$ along $h$.
It is easy to check that $p$ and $s$ are compatible Poisson brackets,
i.e., the Schouten bracket $\[p,s\]=0$.

A pair of compatible Poisson brackets we call a Poisson pencil.

\begin{cor}
\label{cor2.1}
Let $\A_{t,h}$ be a two parameter $\Uh$ invariant quantization
of a commutative algebra $\A$ such that $\A_{t,0}$ is a
one parameter $\Ug$ invariant quantization of $\A$ with
Poisson bracket $s$. 
Then the $\Uh$ invariant quantization $\A_{0,h}$ has a Poisson bracket
$p$ of the form (\ref{rinv}): $p=f-\{\cdot,\cdot\}_r$,  
where  $f$ is an invariant bracket such that $\[f,f\]=-\ff_\A$ and
compatible with $s$, i.e.,
\be{}
\label{cfs}
\[f,s\]=0.
\ee{}
\end{cor}

\begin{proof}
For the two parameter quantization, the Poisson brackets
$p$ and $s$ form a Poisson pencil, hence must be compatible.
Also, $s$ is a $\Ug$ invariant bracket, so that $s$ is compatible
with the $r$-matrix bracket $\{\cdot,\cdot\}_r$. 
It follows from (\ref{rinv}) that $s$ has to be compatible
with $f$.
\end{proof}  

In what follows, we shall often call $\A_{t,h}$ a $\Uh$ invariant 
quantization (or double quantization) of the invariant Poisson bracket $s$,
or of the Poisson pencil $s$ and $p$.

\begin{rem}\label{rem2.4}
As we have seen in Subsection \ref{ss2.3}, to construct a $\Uh$ invariant quantization
of $\A$ is the same that to construct a $\Ug$ invariant $\Phi_h$ associative
quantization of $\A$. We shall see that the last problem
often turns out to be simpler (see Subsection 4.3).
We observe that if $p=f-\{\cdot,\cdot\}_r$ is an admissible 
Poisson bracket for $\Uh$ invariant quantization, then the invariant
bracket $f$ with the property $\[f,f\]=-\ff_\A$ may be
considered as a ``Poisson bracket'' of quantization in the category
with $\Phi_h$ defining the associativity constraint.
Also, the pair $f$, $s$ is a Poisson pencil in that category.
\end{rem}

\section{Double quantization on coadjoint representations}
\label{s3}
In this section we study a two parameter (or double) quantization
on coadjoint representations of simple Lie algebras.

Let $\g$ be a complex Lie algebra. 
Then, the symmetric algebra $S\g$ can be considered as a function algebra
on $\g^*$. The algebra $U(\g)$ is included in the family of algebras
$(S\g)_t=T(\g)[t]/J_t$, where $J_t$ is  the ideal generated by the elements
of the form $x\ot y-\sigma(x\ot y)-t[x,y]$, $x,y\in\g$, $\sigma$ is
the permutation. By the PBW theorem, $(S\g)_t$ is a free module
over $\C[t]$. 
We have $(S\g)_0=S\g$, so this
family of quadratic-linear algebras gives a $U(\g)$ invariant 
quantization of $S\g$ by the Lie bracket $s$.

It turns out that for $\g=sl(n)$ this picture can be extended to 
the quantum case, \cite{Do2}.
Namely, there exist deformations, $\sigma_h$ and $[\cdot,\cdot]_h$, of 
both the mappings
$\sigma$ and $[\cdot,\cdot]$ such that
the two parameter family of algebras
$(S\g)_{t,h}=T(\g)[[h]][t]/J_{t,h}$, where $J_{t,h}$ is the ideal generated 
by the elements
of the form $x\ot y-\sigma_h(x\ot y)-t[x,y]_h$, $x,y\in\g$,
gives a $\Uh$ invariant quantization of the Lie bracket $s$ on $\g^*$.
In this case, the corresponding bracket $f$ from (\ref{cfs}) is a quadratic
bracket which is, up to a factor, a unique nontrivial invariant map 
$\wedge^2\g\to S^2\g$.

We shall show that for other simple Lie algebras, 
double quantizations of the Lie brackets do not exist.

We give two constructions of the algebra
 $(S\g)_{t,h}$.
The first construction uses an idea from the paper \cite{LS} on 
a quantum analog of Lie algebra for $sl(n)$.
The second construction using the so called reflection equations (RE), \cite{KS}, \cite{Maj},
is presented in Remark \ref{rem3.1}..

\subsection{Quantum Lie algebra for $U_h(sl(n))$}
\label{ss3.1}
Let $\Uh$ be a quantized universal enveloping algebra for a Lie algebra $\g$.
We consider $\Uh$ as a $\Uh$ module with respect to the left adjoint action:
$\ad(x)y=x_1y\gamma(x_2)$, where $x,y\in \Uh$, $\De_h(x)=x_1\ot x_2$
(summation implicit).

There were  attemptions to define quantum Lie algebras as
deformed standard classical embeddings of $\g$ into $\Uh$ obeying
some additional properties,
\cite{DG}, \cite{LS}. 

In the classical case, there is probably the following way (not using comultiplication)
to distinguish the standard embedding $\g\to \Ug$ from other invariant embeddings:
with respect to this embedding $\Ug$ is a quadratic-linear algebra.
So, we give the following (working) definition of quantum Lie algebras.

\begin{defn}\label{def3.1}
Let $\g_h$ be a subrepresentation of $\Uh$, which is
a deformation of the standard embedding of $\g$ in $\Ug$.
We call $\g_h$ a quantum Lie algebra, if the kernel of the induced
homomorphism $T(\g_h)\to\Uh$ is defined by (deformed) quadratic-linear relations.
\end{defn}

We are going to show that the quantum Lie algebra exists in case $\g=sl(n)$.
On the other hand, if such an algebra exists for some Lie algebra $\g$,
then a double quantization of the Lie bracket on $\g^*$ also exists.
But, as we shall see, no double quantization exists for simple $\g\neq sl(n)$.
So, among simple finite dimensional Lie algebras, only $sl(n)$ has a quantum Lie algebra
in our sense.

Our construction is the following.
Let $R=R^\prime_i\ot R^{\prime\prime}_i\in \Uh\ot \Uh$ 
(completed tensor product) be the R-matrix 
(summation by $i$ is assumed).
It satisfies the properties \cite{Dr2}
\be{} \label{Delo}
\Delta_h^\prime(x)=R\Delta_h(x)R^{-1}, \quad x\in \Uh,
\ee{}
where $\Delta_h$ is the comultiplication in $\Uh$ and $\Delta_h^\prime$ is  
the opposite one,
\be{} \label{six}
(\Delta_h\ot 1)R&=&R^{13}R^{23}=R^\prime_i\ot R^\prime_j\ot R^{\prime\prime}_i
R^{\prime\prime}_j \notag\\
(1\ot\Delta_h)R&=&R^{13}R^{12}=R^{\prime}_iR^{\prime}_j\ot R^{\prime\prime}_j
\ot R^{\prime\prime}_i, 
\ee{}
and
\be{} \label{coun}
(1\ot \varepsilon)R=(\varepsilon\ot 1)R=1\ot 1, 
\ee{}
where $\varepsilon$ is the counit in $\Uh$.

Consider the element $Q=Q_i^\prime\ot Q_i^{\prime\prime}=R^{21}R$. 
It follows from (\ref{Delo}) that $Q$
commutes with elements from $\Uh\ot \Uh$ of the form $\Delta_h(x)$. This is
equivalent for $Q$ to be invariant under the adjoint action of $\Uh$
on  $\Uh\ot \Uh$.

Let $V$ be an irreducible finite dimensional representation of $\Uh$ and 
$\rho:\Uh\to\End(V)$ the corresponding map of algebras.
Consider the dual space $\End(V)^*$ as a left $\Uh$ module setting
$$(x\varphi)(a)=\varphi(\gamma(x_{(1)})ax_{(2)}), $$
where $\varphi\in\End(V)^*$, $a\in\End(V)$, $\Delta_h(x)=x_{(1)}\ot x_{(2)}$
in Sweedler notions, and $\gamma$ denotes the antipode in $\Uh$.

Consider the map 
\be{}\label{mapE}
f:\End(V)^*\to\Uh 
\ee{}
defined as
$\varphi\mapsto \varphi(\rho(Q^{\prime}_i)Q^{\prime\prime}_i$. 
From the invariance of
$Q$ it follows that $f$ is a $\Uh$ equivariant map, so $\bL=\Im(f)$ is a
$\Uh$ submodule. 

It follows from (\ref{six}) that $\bL$ is a left coideal in $\Uh$, i.e.,
$\Delta(x)\in \Uh\ot \bL$ for any $x\in \bL$. 
Indeed, $Q=R^{\prime\prime}_iR^{\prime}_j\ot R^{\prime}_iR^{\prime\prime}_j$.
Applying  
(\ref{six}) we obtain
\be{}\label{DeQ}
(1\ot\Delta_h)R^{21}R=R^{\prime\prime}_i
R^{\prime\prime}_jR^{\prime}_kR^{\prime}_l\ot R^{\prime}_iR^{\prime\prime}_l\ot
R^{\prime}_jR^{\prime\prime}_k
\ee{}
Let $\varphi\in \End(V)^*$. Define $\psi_{il}\in\End(V)^*$ setting
$\psi_{il}(a)=\varphi(R^{\prime\prime}_iaR^{\prime}_l)$ for $a\in\End(V)$.
Then
$\Delta\varphi(R^{\prime\prime}_iR^{\prime}_j)R^{\prime}_iR^{\prime\prime}_j=
R^{\prime}_iR^{\prime\prime}_l\ot \psi_{il}(R^{\prime\prime}_j
R^{\prime}_k)R^{\prime}_jR^{\prime\prime}_k,$
which obviously belongs to $\Uh\ot\bL$. 
\smallskip

Recall, \cite{Dr2}, that $R=F^{21}_he^{\frac{h}{2}\tb}F^{-1}_h$. 
Here $\tb=\sum_i t_i\ot t_i$ is the split Casimir, where 
$t_i$ form an orthonormal
basis in $\g$ with respect to the Killing form, 
$F=1\ot 1+\frac{h}{2}r+o(h)$ (see (\ref{F})),
and $r$ is a classical  $r$-matrix. 
Therefore, 
\be{} \label{rmt}
Q=R^{21}R=Fe^{h\tb}F^{-1}=1\ot 1+h\tb+\frac{h^2}{2}(\tb^2+[r,\tb])+o(h^2).
\ee{}

Denote by $\Tr$ the unique (up to a factor) invariant element in $\End(V)^*$.
Let $Z_0=\rho_0(\g)$, and denote by $Z_h$ some $\Uh$ invariant deformation of
$Z_0$ in $\End(V)$.  
Then we have a decomposition $\End(V)=I\oplus Z_h\oplus W$, where $I$ is 
the one dimensional invariant subspace generated by the identity map,
$W$ is a complement to $I\oplus Z_h$ invariant subspace. 
This gives a decomposition $\End(V)^*=I^*\oplus Z_h^*\oplus W^*$ where
$W^*$ consists of all the elements which are equal to zero on $I\oplus Z_h$.  
The space $I^*$ is generated by $\Tr$, and after normalizing in such a way
that $\Tr(\id)=1$, we obtain that $C_V=f(\Tr)$ is of the form
\be{} \label{rcv}
C_V=\Tr(\rho(Q^\prime_i))Q^{\prime\prime}_i=1+h^2c+o(h^2),
\ee{}
where $c$ is an invariant element of $U(\g)$. It follows from (\ref{coun})
that $\varepsilon(C)=1$.

From (\ref{rmt}) follows that the elements of $f(Z_h^*)$ have the form
\be{}
z=hx+o(h), \quad x\in \g,
\ee{}
hence the subspace $L_1=h^{-1}f(Z_h^*)$ forms a subrepresentation of $\Uh$
with respect to the left adjoint action of $\Uh$ on itself, which is a deformation
of the standard embedding of $\g$ into $U(\g)$. It follows from (\ref{coun})
that $\varepsilon(L_1)=0$.

The elements from $f(W^*)$ have the form $w=h^2b+o(h^2)$
and $\varepsilon(W^*)=0$. 
Denote $L_2=h^{-2}f(W^*)$.

So,  
$\bL=\C C_V\oplus hL_1\oplus h^2L_2=\C C_V+hL$,
where $L=L_1\oplus hL_2$.
Since $\bL$ is a left coideal in $\Uh$, for any $x\in \bL$ we have
\be{*} 
\Delta_h(x)=x_{(1)}\ot x_{(2)}=z\ot C_V+v\ot x^{\prime},
\ee{*}
where $z,v\in \Uh$, $x^{\prime}\in L$.
Applying to the both hand sides $(1\ot\varepsilon)$ and multiplying
we obtain $x=x_{(1)}\vep(x_{(2)})=z\vep(C_V)+v\vep(x^{\prime})=z$. So, $z$ has
to be equal to $x$. and we obtain
\be{} \label{qlr1}
\Delta_h(x)=x_{(1)}\ot x_{(2)}=x\ot C_V+v\ot x^{\prime}, 
\quad x,x^{\prime}\in L.
\ee{}
From (\ref{qlr1}) we have for any $y\in L$ 
\be{} \label{qlr2}
xy=x_{(1)}y\gamma(x_{(2)})x_{(3)}=x_{(1)}y\gamma(x_{(2)})C_V+v_{(1)}y
\gamma(v_{(2)})x^{\prime}.
\ee{}

Introduce the following maps:
\be{} \label{qlr'}
\sigma_h^{\prime}:&L\ot L\to L\ot L,&\quad x\ot y\mapsto 
v_{(1)}y\gamma(v_{(2)})\ot x^{\prime},  \notag \\  
\ [\cdot,\cdot]_h^{\prime}:&L\ot L\to L, &\quad x\ot y\mapsto x_{(1)}y
\gamma(x_{(2)}).
\ee{}
We may rewrite (\ref{qlr2}) in the form
\be{} \label{qlr3}
m(x\ot y-\sigma_h^{\prime}(x\ot y))-[x,y]_h^{\prime}C_V=0.
\ee{}

Observe now that, as follows from (\ref{rcv}),  
$C_V$ is an invertible element in $\Uh$. Put $P=C^{-1}_V$.
Transfer the maps (\ref{qlr'}) to the space $P\cdot L$, i.e., define
\be{*}
\sigma_h(Px,Py)&=&(P\ot P)\sigma^{\prime}_h(x,y), \\
{[}Px,Py]_{h} &=&P[x,y]^{\prime}_h.
\ee{*}
From (\ref{qlr1}) we obtain
\be{}
P_{(1)}x_{(1)}\ot P_{(2)}x_{(2)})=P_{(1)}x\ot P_{(2)}C_V+P_{(1)}v\ot 
P_{(2)}x^{\prime}.
\ee{}
Using this relation and taking into account that $P$ commutes with all 
elements from $\Uh$, we obtain as in (\ref{qlr2}) 
\be{}\label{PxPy}
PxPy&=&P_{(1)}x_{(1)}Py\gamma(x_{(2)})\gamma(P_{(2)})P_{(3)}x_{(3)}= \\
& &P_{(1)}x_{(1)}Py\gamma(x_{(2)})\gamma(P_{(2)})P_{(3)}C_V+P_{(1)}v_{(1)}Py 
\gamma(v_{(2)})\gamma(P_{(2)})P_{(3)}x^{\prime}= \\
& & P[x,y]^{\prime}_h+P^2m\sigma_h^{\prime}(x\ot y)= 
 [Px,Py]_h+m\sigma_h(Px\ot Py)). 
\ee{}
This equality may be written as
\be{} \label{qlr4}
m(x\ot y-\sigma_h(x\ot y))-[x,y]_h=0, \quad x,y\in C^{-1}_VL. 
\ee{}

Define $L_V=C^{-1}_VL$.
Let $T(L_V)=\oplus_{k=0}^\infty L_V^{\ot k}$ be the tensor algebra over $L_V$.
Notice, that $T(L_V)$ is not supposed to be completed in $h$-adic topology.
Let $J$ be the ideal in $T(L_V)$ generated by the relations
\be{} \label{rel0}
(x\ot y-\sigma_h(x\ot y))-[x,y]_h, \quad x,y\in L_V.
\ee{}
Due to (\ref{qlr4}) we have a homomorphism of algebras over $\C[[h]]$
\be{}
\label{funm}
\psi_h:T(L_V)/J\to \Uh, 
\ee{}
extending the natural embedding $L_V\to \Uh$ of $\Uh$ modules..

Now we can prove

\begin{propn}
\label{prop3.1}
For $\g=sl(n)$ the quantum Lie algebra exists.
\end{propn}

\begin{proof}

Apply the above construction to $V=\C^n[[h]]$,
the deformed basic representation of $\g$. In this case $\End(V)=I\oplus Z_h$, 
where $Z_h$ is a deformed adjoint representation.
So, $\g_h=L_V=h^{-1}C_V^{-1}f(Z^*_h)$ is a deformation of the standard 
embedding of $\g$ in $U(\g)$. It is easy to see that in this case
$\sigma_h$ is a deformation of the usual permutation: 
$\sigma_0(x\ot y)=y\ot x$, and $[\cdot,\cdot]_h$ is a deformation of the
Lie bracket on $\g$: $[x,y]_0=[x,y]$, $x,y\in \g\subset U(\g)$.
Hence, at $h=0$, the quadratic-linear relations (\ref{rel0}) are exactly the defining relations for
$\Ug$, therefore the map (\ref{funm}) is an isomorphism at $h=0$.
It follows that (\ref{funm}) is an embedding. (Actually, (\ref{funm}) is essentially
an isomorphism, i.e., it is an isomorphism after completion of $T(L_V)$
in $h$-adic topology.) So, the kernel of the map $T(L_h)\to \Uh$  is defined
by the quadratic-linear relations (\ref{rel0}).
\end{proof}

\begin{rem}\label{remRE}
Quadratic-linear relations (\ref{rel0})
can be obtained in another way.
Note that equation (\ref{DeQ}) may be rewritten as
\be{}\label{DeQr}
(1\ot\De_h)Q=R_{21}Q_{13}R_{12}.
\ee{}
Since $Q$ commutes with all elements of the form $\De_h(x)$,
$x\in \Uh$, one derives from (\ref{DeQr}):
\be{}\label{pRE}
Q_{23}R_{21}Q_{13}R_{12}=R_{21}Q_{13}R_{12}Q_{23}.
\ee{}
Consider the element  $Q_\rho=\rho(Q_1)\ot Q_2$ as
a $\dim(V)\times\dim(V)$ matrix with the entries from $\Uh$.
Applying to (\ref{pRE}) operator $\rho\ot \rho\ot 1$,
we obtain the following relation for $Q_\rho$:
\be{}\label{pREr}
(Q_\rho)_2\overline{R}_{21}(Q_\rho)_1\overline{R}=
\overline{R}_{21}(Q_\rho)_1\overline{R}(Q_\rho)_2,
\ee{}
where $\overline{R}=(\rho\ot\rho)R$ is a number matrix, 
the Yang-Baxter operator in $V\ot V$.
Replacing in this equation $\overline{R}$ by $S=\sigma\overline{R}$,
we obtain that  the matrix $Q_\rho$ satisfies the following
reflection equation (RE):
\be{}\label{RE}
(Q_\rho)_2S(Q_\rho)_2S=S(Q_\rho)_2S(Q_\rho)_2.
\ee{}

It is clear that the entries of the matrix $Q_\rho$ generate
the image of the map (\ref{mapE}).
From (\ref{rcv}) follows that $Q_\rho$ has the form
\be{}\label{v1}
Q_\rho=\Id_V\cA_V+hB^\prime,
\ee{}
where $B^\prime$ has the form $B^\prime=\sum D_i\ot b_i$,
$D_i$ belong to the complement to $\C\Id_V$ submodule in
$\End(V)$ and $b_i\in\Uh$. 
Note that the entries of the matrix $B^\prime$ form the subspace $L$,
whereas the entries of $B=\cA_V^{-1}B^\prime$ form
the subspace $L_V$ from (\ref{qlr4}).
From (\ref{v1}) we obtain
\be{}\label{v2}
\cA_V^{-1}Q_\rho=\Id+hB.
\ee{}
Since the element $\cA_V^{-1}$ belongs to the center of $\Uh$,
the matrix $\cA_V^{-1}Q_\rho$ obeys the RE (\ref{RE}) as well. 
So, $B$ satisfies the relation
\be{}\label{REm}
(\Id+hB)_2S(\Id+hB)_2S=S(\Id+hB)_2S(\Id+hB)_2. 
\ee{}
One checks that (\ref{REm}), considered as a qudratic-linear relations
for indetermined entries of $B$, is equivalent to (\ref{rel0}) in the case $\g=sl(n)$.
\end{rem}

\subsection{Double quantization on $sl(n)^*$}

Introduce a new variable, $t$, and consider a homomorphism of algebras,
$T(L_V)[t]\to \Uh[t]$, which extends the embedding $t\cdot\imath:L_V[t]\to \Uh[t]$,
where $\imath$ stands for the standard embedding $L_V\to \Uh$.
From (\ref{qlr4}) follows that $t\cdot\imath$ factors through the homomorphism
of algebras over $\C[[h]][t]$
\be{} \label{phit}
\phi_{t,h}:T(L_V)[t]/J_t\to \Uh[t],
\ee{} 
where $J_t$ is the ideal generated by the relations
\be{} \label{relt}
(x\ot y-\sigma_h(x\ot y))-t[x,y]_h, \quad x,y\in L_V.
\ee{}

\begin{propn}\label{prop3.2}
For $\g=sl(n)$ the algebra $(S\g)_{t,h}=T(L_V)[t]/J_t$ is a double quantization of 
the Lie bracket on $S\g$.
\end{propn}

\begin{proof}
Since in this case $L_V=\g_h$, from Proposition \ref{prop3.1} follows that 
(\ref{phit}) is a monomorphism at $t=1$. 
Due to the PBW theorem the algebra $\Im(\phi_{t,h})$ 
at the point $h=0$ is a free $\C[t]$-module and is equal to
\be{} \label{qst}
(S\g)_t=T(\g)/\{x\ot y-y\ot x-t[x,y]\}. 
\ee{}
For $t=0$ this algebra is the
symmetric algebra $S\g$, the algebra of algebraic functions on $\g^*$. 
For $t\neq 0$ this algebra 
is isomorphic to $U(\g)$. Since $\Uh$ is a free $\C[[h]]$-module,
it follows 
that $\phi_{t,h}$ in (\ref{phit}) is a monomorphism of algebras over 
$\C[[h]][t]$
and $\Im(\phi_{t,h})$ is a free $\C[[h]][t]$-module isomorphic to
\be{} \label{qsth}
(S\g)_{t,h}=T(\g_h)[t]/\{x\ot y-\sigma_h(x\ot y)-t[x,y]_h\}. 
\ee{}
It is clear that $(S\g)_t=(S\g)_{t,0}$ is the standard quantization of the
Lie bracket on $\g^*$.
\end{proof}

Call the algebra 
\be{} \label{qsh}
(S\g)_h=(S\g)_{0,h}=T(\g_h)/\{x\ot y-\sigma_h(x\ot y)\}
\ee{}
a quantum symmetric algebra (or quantum polynomial algebra on $\g^*$). 
It is a free $\C[[h]]$ module and a quadratic algebra equal to $S\g$ at $h=0$.

\begin{rem}\label{rem3.2} Up to now, all our constructions were considered for
the quantum group in sense of Drinfeld, $\Uh$, defined over $\C[[h]]$.
But one can deduce the results above for
the quantum group in sense of Lusztig, $U_q(\g)$, defined over
the algebra $\C[q,q^{-1}]$.
We show, for example, how to obtain the quantum
symmetric algebra over $\g$.
Let $E$ be a Grassmannian consisting of subspaces  
in $\g\ot\g$ of dimension equal to $\dim(\wedge^2\g)$,
and $\Z$ the closed algebraic subset of $E$ consisting of 
subspaces $J$ such that
$\dim(E\ot J\cap J\ot E)\geq\dim(\wedge^3\g)$.
Let $\X$ be 
the algebraic subset in $\Z\times(\C\setminus 0)$   
consisting of points $(J,q)$ such that $J$
is invariant under the action of $U_q(\g)$. The projection
$\pi:\X\to\C\setminus 0$ is a proper map.
It is clear that the fiber of this projection over $q=1$ contains 
the point corresponding
to the symmetric algebra $S\g$ as an isolated point,
because 
there are no quadratic $U(\g)$ invariant Poisson
brackets on $S\g$.

As follows from the existence of $(S\g)_h$ (completed situation at $q=1$),
the dimension of $\X$ is equal to $1$. Hence, the projection 
$\pi:\X \to\C\setminus 0$ is a covering. For $x\in\X$ let $J_x$ be the
corresponding subspace in $\g\ot\g$ and $(S\g)_x=T(\g)/\{J_x\}$
the corresponding quadratic algebra. Due to the projection $\pi$, the
family $(S\g)_x$, $x\in\X$, is a module over $\C[q,q^{-1}]$. Since
$J_x$ is $U_{p(x)}(\g)$ invariant, $(S\g)_x$ is a $U_{p(x)}(\g)$ invariant 
algebra.
Hence, after possible deleting from
$\X$ some countable
set of points, we obtain a  family of quadratic algebras with the same
dimensions of graded components as $S\g$. So, the family $(S\g)_x$, $x\in \X$ 
can be considered as a quantum symmetric
algebra over $U_q(\g)$.

Note also that the family $(S\g)_h$ can be thought of as a 
formal section of the map $\pi:\X\to (\C\setminus 0)$ 
over the formal neighborhood of point $q=1$.
It follows that there is also
an analytic section of $\pi$ over some neighborhood, $U$, of the point $q=1$.
If $(S\g)_h$ is a quantization with Poisson bracket $f-\{\cdot,\cdot\}_r$
(see Proposition \ref{prop2.2}), then 
a quantization with Poisson bracket $-f-\{\cdot,\cdot\}_r$ 
gives another section of $\pi$ over $U$.
Hence, in a neighborhood of the ``classical'' point $x_0\in\X$, $\pi(x_0)=1$,
the space $\X$ has a singularity of type ``cross''.
 
\end{rem}

\subsection{Poisson pencil corresponding to $(S\g)_{t,h}$}

Let $\g=sl(n)$ and $(S\g)_{t,h}$  be the double quantization from
Proposition \ref{prop3.2}.

\begin{propn}\label{prop3.3}
The Poisson pencil corresponding to the quantization $(S\g)_{t,h}$ consists of
two compatible Poisson brackets: 

$s$ (along $t$) is the Lie bracket;

$p$ (along $h$) is a quadratic Poisson bracket of the form
$p=f-\{\cdot,\cdot\}_r$, where $f$ is an invariant quadratic bracket
which is a unique up to a factor
invariant map $f:\wedge^2\g\to S^2\g$, and $\{\cdot,\cdot\}_r$
is the $r$-matrix bracket. Moreover, $\[s,f\]=0$ and
$ \[f,f\]=-\overline{\ff}$, where $\overline{\ff}$
has the form $\overline{\ff}(a,b,c)=[\ff_1,a][\ff_2,b][\ff_3,c]$,
and $\ff=\ff_1\wedge\ff_2\wedge\ff_3=\[r,r\]$. 
Recall that $\ff$ is a unique up to a factor invariant element of $\wedge^3\g$.

\end{propn}

\begin{proof}

That $s$ coincides with the Lie bracket is obvious from (\ref{qst}).
From Corollary~\ref{cor2.1} we have  $p=f-\{\cdot,\cdot\}_r$. 
Since $(S\g)_h$ is a quadratic algebra over $\C[[h]]$, $p$ must be a quadratic bracket.
But the $r$-matrix bracket $\{\cdot,\cdot\}_r$ is quadratic, too. Hence,
$f$ must be a quadratic invariant bracket. There is only one possibility
for such a bracket: it must be a unique (up to a factor) nontrivial
invariant map $f:\wedge^2\g\to S^2\g$. Now apply
Proposition \ref{prop2.2} and Corollary \ref{cor2.1}.
\end{proof}

Consider now the quadratic bracket $f$ in more detail.

We say that a $k$-vector field, $g$, on a manifold  $M$ is strongly restricted
on a submanifold $N\subset M$ if at any point of $N$ the polyvector $g$
can be presented as an exterior power of tangent vectors 
to $N$.

Consider the coadjoint action of the Lie group $G=SL(n)$ on $\g^*=sl(n)^*$.
We want to prove that the bracket $f$ is strongly restricted on
any orbit of $G$ in $\g$. 
It turns out that there is the following general fact.

\begin{propn}\label{prop3.4}
Let $G$ be a semisimple Lie group with its Lie algebra $\g$, 
$s=[\cdot,\cdot]$  the Lie bracket on $\g^*$.
Let $f=\{\cdot,\cdot\}$ be an invariant bracket on $\g^*$ such that the Schouten
bracket $\[s,f\]$ is a three-vector field, $\psi$, strongly restricted on an orbit $\O$
of $G$. 
Then $f$ is strongly restricted on $\O$.
\end{propn}

\begin{proof}
Let $x,y,z\in\g$.
The invariance condition for $\{\cdot,\cdot\}$ means:
\be{}\label{invbrcond}
[x,\{y,z\}]=\{[x,y],z\}+\{y,[x,z]\}.
\ee{}
The Schouten bracket $\[s,f\]$ is:
\be{*}\label{schoutcond}
[x,\{y,z\}]+[y,\{z,x\}]+[z,\{x,y\}]+\\
\{x,[y,z]\}+\{y,[z,x]\}+\{z,[x,y]\}=\psi(x,y,z).
\ee{*}
In the left hand side of this expression, the 1-st, 5-th, and 6-th terms are canceled
due to (\ref{invbrcond}), and we have
\be{*}
[y,\{z,x\}]+[z,\{x,y\}]+\{x,[y,z]\}=\psi(x,y,z).
\ee{*}
Putting in this equation instead of  $[y,\{z,x\}]$ its expression from 
(\ref{invbrcond}), i.e., $\{[y,z],x\}+\{z,[y,x]\}$, we obtain,
since the term $\{x,[y,z]\}$ is canceled:
\be{}\label{scomp}
\{z,[x,y]\}=[z,\{x,y\}]+\psi(x,y,z).
\ee{}
Now observe that, due to the Leibniz rule, equation (\ref{scomp}) is valid
for any $z\in S\g$. To prove the proposition, it is sufficient
to show that if $z$ belongs to the ideal $I_\O$ defining the orbit $\O$,
then $\{z,u\}$ also belongs to this ideal. Again, due to the Leibniz rule,
it is sufficient to show this for $u\in\g$.
Since $\g$ is semisimple, there are elements $x,y\in\g$ such that
$[x,y]=u$. We have from (\ref{scomp})
\be{*}
\{z,u\}=\{z,[x,y]\}=[z,\{x,y\}]+\psi(x,y,z).
\ee{*}
But $[z,\{x,y\}]\in I_\O$, since the Lie bracket is restricted
on any orbit, $\psi(x,y,z)\in I_\O$ by hypothesis of the proposition.
So, $\{z,u\}\in I_\O$.
\end{proof}

As a consequence we obtain
\begin{propn}\label{prop3.5}
Let $\g=sl(n)$. Then the bracket $f$ from Proposition \ref{prop3.3} 
is strongly restricted on any orbit of $SL(n)$.
\end{propn}

\begin{proof}
Follows from Propositions \ref{prop3.3} and \ref{prop3.4}.
\end{proof}

\begin{rem}\label{remPL1}
According to Remark \ref{remPL}, this Proposition shows
that in case $G=SL(n)$ any orbit in coadjoint representation
has a Poisson bracket $p=f-r_M$ such that the pair $(M,p)$
becomes a $(G,\tilde{r})$-Poisson manifold.
\end{rem} 

\begin{rem}\label{rem3.1}
Recall that in case $\g=sl(n)$ the tensor square $\g\ot\g$, considered as a
representation of $\g$, has a decomposition into irreducible components
which are contained in  $\g\ot\g$ with multiplicity one, except of the component
isomorphic to $\g$ having multiplicity two. Moreover, both
the symmetric and skew-symmetric
parts of $\g\ot\g$ contain  components, $\g^1$ and $\g^2$, 
isomorphic to $\g$. Hence, the bracket $f$ takes $\g^2$ onto $\g^1$
and all the other components to zero.

For $\g$ simple not equal to $sl(n)$, the decomposition of  $\g\ot\g$ 
is multiplicity free, hence non-trivial invariant maps 
$\wedge^2\g\to S^2\g$ do not exist at all. It follows that for $\g\neq sl(n)$,
there do not exist quadratic algebras $(S\g)_{h}$ which are
$\Uh$ invariant quantizations of $S\g$.
\end{rem}
\begin{conjecture}
Prove that there exist no one parameter $\Uh$ invariant quantizations of $S\g$
(not necessarily in the class of quadratic algebras) for all
simple Lie algebras $\g\neq sl(n)$.
\end{conjecture}

Now we prove that for simple $\g\neq sl(n)$, the double quantization 
does not exist (not necessarily in the class of quadratic-linear algebras).

\begin{propn}\label{prop3.6}

Let $\g$ be a simple finite dimensional Lie algebra not
equal to $sl(n)$. Then a $\Uh$ invariant quantization of the Lie
bracket on $\g^*$ does not exist.
\end{propn}

\begin{proof}
If such a quantization exists, then from Corollary \ref{cor2.1} follows
that there exists an invariant bracket $f$ on $\g^*$ such that
$\[s,f\]=0$ and $\[f,f\]=-\overline{\ff}$. Here $s$ is the Lie bracket and
$\overline{\ff}$ is the three-vector field induced by $\ff$ (see Proposition
\ref{prop3.3}). We show that such  $f$ does not exist.
Observe that $\overline{\ff}$ has type $(3,3)$, i.e.,
is a sum of terms of the view $b\,\partial_x\wedge\partial_y\wedge\partial_z$,
where $b$ is a homogeneous polynomial of degree $3$.
Observe also that the Schouten bracket of two polyvector fields
of degrees $(i,j)$ and $(k,l)$ is a polyvector field of degree $(i+k-1,j+l-1)$.
We shall write $i$ for degree $(i,j)$ when the second number, $j$, is
clear from context.

It is obvious that on $\g$ there are no invariant bivector fields of degree $0$
and,  up to a factor, there is a unique invariant bivector field
of degree $1$, the Lie bracket $s$ itself.
Since $\g\neq sl(n)$, there are no bivector fields of degree $2$ 
(see Remark \ref{rem3.1}). Therefore, $f$ must be of the form:
$f=s+f_1$, where $f_1$ is a bracket of degree $\geq 3$.
Since $f$ is compatible with $s$ and $\[f,f\]=-\overline{\ff}$, it must be
$\[ f_1,f_1\]=-\overline{\ff}$. But it is impossible, because
$\[ f_1,f_1\]$ has at least degree $5$. 
\end{proof}

\subsection{Quantum de Rham complex on $(sl(n))^*$}

Consider the algebra $\Omega^\bullet$ of differencial forms on $\g^*$ with polynomial
coefficients. This is a graded differential algebra with  differential $d$ of
degree $1$
which forms the de Rham complex
\be{}\label{cdR}
\ee{}
where $\Omega^k$ is the space of $k$-forms with polynomial
coefficients.

We call a complex over $\C[[h]]$
\be{}\label{qdR}
\Omega_h^\bullet: (S\g)_h \tor{d_h} \Omega_h^1 \tor{d_h} \Omega_h^2\tor{d_h}\cdots
\ee{}
a quantum (deformed) de Rham complex if it consists of $\Uh$ invariant
topologically
free modules over $\C[[h]]$
and coincides with (\ref{cdR}) at $h=0$.

\begin{propn}\label{prop3.7}
Let $\g=sl(n)$. 
Then the quantized polynomial algebra $(S\g)_h$ from (\ref{qsh})
can be included in a $\Uh$ invariant  graded differential algebra,
$\Omega_h^\bullet$, which form a 
quantum de Rham complex (\ref{qdR}).
\end{propn}

\begin{proof}
First of all, define a quantum exterior algebra, $(\Lambda\g)_h$, an
algebra of differential forms with constant coefficients.
Let us modify the operator $\sigma_h$ from (\ref{qsh}). 
Since the representation $\g^*_h$ is 
isomorphic to $\g_h$, there exists a $\Uh$ invariant bilinear form on $\g_h$,
deformed Killing form. 
This form can be naturally extended to all tensor degrees $\g_h^{\ot k}$.
Let $W_h^2$ be the $\C[[h]]$ submodule in $\g_h\ot \g_h$
orthogonal to $V_h^2=\Im(\id\ot\id-\sigma_h)$. 
Define an operator $\bsig$ on $\g_h\ot\g_h$ in such a way that it has
the eigenvalues $-1$ on $V_h^2$ and
$1$ on $W_h^2$. It is clear  that $V_h^2$ and $W_h^2$ are 
deformed skew symmetric
and symmetric subspaces of $\g\ot \g$. 

%%%%
Now observe that the third graded component in the quadratic algebra
$(S\g)_h$ is the quotient of $\g_h^{\ot 3}$ by the submodule 
$V_h^2\ot\g_h+\g_h\ot V_h^2$, hence this submodule and, 
therefore, the submodule
$V_h^2\ot\g_h\cap \g_h\ot V_h^2$ are direct submodules in $\g_h^{\ot 3}$, i.e., 
they have
complement submodules. As the complement submodules one can
choose the submodules $W_h^2\ot\g_h\cap \g_h\ot W_h^2$ and 
$W_h^2\ot\g_h+\g_h\ot W_h^2$, respectively, since they are complement 
at the point $h=0$
and $W_h^2$ is orthogonal to $V_h^2$ with respect to the Killing form extended
to $\g_h\ot\g_h$. 
Hence, $W_h^2\ot\g_h+\g_h\ot W_h^2$ is a direct submodule. 
Also, the symmetric algebra $S\g$ is Koszul. From
a result of Drinfeld, \cite{Dr3} (see also \cite{DM}),  follows that the quadratic algebra 
$(\Lambda\g)_h=T(\g_h)/\{W_h^2\}$ is a free $\C[[h]]$ module, i.e., is
a $\Uh$-invariant deformation of the exterior algebra $\Lambda\g$.

Call $(\Lambda\g)_h$ a quantum exterior algebra over $\g$.

Define a quantum algebra of differential forms over $\g^*$ as the tensor
product $\Omega_h^\bullet=(S\g)_h\ot(\Lambda\g)_h$ in the tensor category of
representations of the quantum group $\Uh$. The multiplication of two
elements $a\ot\alpha$ and $b\ot\beta$ looks like $ab_1\ot\alpha_1\beta$,
where $b_1\ot\alpha_1=S(\alpha\ot b)$ for $S=\sigma R$ being the permutation
in that category. So, $\Omega_h^k=(S\g)_h\ot (\Lambda^k\g)_h$. 

As in the classical case, the algebras $(S\g)_h$ and $(\Lambda\g)_h$
can be embedded in $T(\g_h)$ as a graded submodules
in the following way. Call the submodule 
$W^k_h=(W_h^2\ot\g_h\ot\cdots\ot\g_h)\cap
(\g_h\ot W_h^2\ot\g_h\ot\cdots\ot\g_h)\cap\cdots\cap
(\g_h\ot\g_h\ot\cdots\ot W_h^2)$
of $T^k(\g_h)$ a $k$-th symmetric part of $T(\g_h)$.
It is clear that the natural map $\pi_W:T(\g_h)\to (S\g)_h$ restricted
to $W^k_h$ is a bijection onto the $k$-degree component $(S^k\g)_h$ 
of $(S\g)_h$. Denote by $\pi_W^{\prime}:(S^k\g)_h\to W^k_h$
the inverse bijection. Similarly we define $V^k_h$, the $k$-th 
skew symmetric part of $T(\g_h)$, and the bijection
$\pi_V^{\prime}:(\Lambda^k\g)_h\to V^k_h$.

Now, define a differential $d_h$ in $\Omega_h^\bullet$ as a homogeneous
operator of degree $(-1,1)$. It acts on an element, $a\ot\omega$,
of degree $(k,m)$ in the following way. Let $a\ot\omega=
(a_1\ot\cdots\ot a_k)\ot(\omega_1\ot\cdots\ot\omega_m)$ be
its realization as an element from $W^k_h\ot V^m_h$.
Then the formula 
\be{}
d_h(a\ot\omega)=(a_1\ot\cdots\ot a_{k-1}\ot\pi_V^\prime\pi_V
(a_k\ot\omega_1\ot\cdots\ot\omega_m)
\ee{}
presents the element $d_h(a\ot\omega)$ through its realization in
$W^{k-1}_h\ot V^{m+1}_h$. It is obvious that $d_h^2=0$.

So, the  graded differential algebra
$\Omega_h^\bullet$ is constructed.
It is easy to see that at the point $h=0$ this algebra coincides with
$\Omega^\bullet$.
\end{proof}

Note that the quantum de Rham complex is exact, because 
it is exact at $h=0$.

\subsection{Restriction of $(S\g)_{t,h}$ on orbits}
\label{ss3.5}
In this section $G=SL(n)$, $\g=sl(n)$.

Let $M$ be an invariant closed algebraic subset in $\g^*$ and $A$ the algebra
of algebraic functions on $M$. The algebra $A$  can be 
presented as a quotient of $S\g$ by some 
ideal, $S\g\to A\to 0$.

We say that the quantization $(S\g)_{t,h}$ can be restricted on
$M$ if there exists a $\Uh$ invariant quantization, $A_{t,h}$, of $A$, which
can be presented as a quotient of $(S\g)_{t,h}$ by some 
ideal, $(S\g)_{t,h}\to A_{t,h}\to 0$.

Note that, on the infinitesimal level, there  are no obstructions 
for $(S\g)_{t,h}$ to be restricted on $M$.
Indeed, the Lie bracket on $\g^*$ is strongly restricted on
any orbit of $G$ and induces the Kirillov-Kostant-Souriau
bracket on $M$. Also , by Proposition \ref{prop3.5},
the bracket $f$ involved in the quantization
along $h$ is also strongly restricted on any orbit.

From \cite{DS1}, one can derive that the problem of restriction 
of $(S\g)_{t,h}$ is solved positively in case $M$ is 
a minimal semisimple orbit, i.e., $M$ is a hermitian
symmetric space. 
 
We are going to show here that the problem also has a positive solution
for $M$ being a maximal semisimple orbit, i.e., if $M$ can be defined as a set of zeros
of invariant functions from $S\g$. Such orbits are the orbits of diagonal
matrices with distinct elements on diagonal.

\begin{propn}\label{prop3.8}
Let $\g=sl(n)$. Then the family $(S\g)_{t,h}$ can be restricted
on any maximal semisimple orbit in $\g^*$.
\end{propn}

\begin{proof}
There exists an isomorphism
of $\Uh$ modules $(S\g)_{h}\to W_h$, where $W_h=\oplus_kW_h^k$, the direct sum 
of the $k$-th symmetric parts of $T(\g_h)$ (see previous Subsection).
Consider the composition $W_h[t]\to T(\g_h)[t]\to(S\g)_{t,h}$, where the last
map appears from (\ref{qsth}). It is
an isomorphism, since it is an isomorphism at the point $h=0$.
It follows that $(S\g)_{t,h}$ is isomorphic to $W_h[t]$ as a $\Uh$-module,

Denote by $\I_{t,h}$ the submodule
of $\Uh$ invariant elements in $(S\g)_{t,h}$. It is obvious that $\I_{t,h}$
is isomorphic to $\oplus_k\I^k_h[t]$, where $\I^k_h$ is the invariant submodule
in $W^k_h$. Hence, $\I_{t,h}$ is a direct free $\C[[h]][t]$ submodule
in $(S\g)_{t,h}$. Moreover, $\I_{t,h}$ is a central subalgebra in 
$(S\g)_{t,h}$. Indeed, for a generic $t$ the algebra $(S\g)_{t,h}$ can be 
invariantly embedded in $\Uh$. But $\ad(\Uh)$ invariant elements in $\Uh$
form the center of $\Uh$. Also, $\I_{t,h}$ as an algebra
is isomorphic to $\I[[h]][t]$ with the trivial action of $\Uh$, where
$\I=\I_{0,0}$, the algebra of invariant elements in $S\g$. This follows from the fact
that $\I$ is a polynomial algebra, \cite{Dix}, and, therefore, admits no  
nontrivial commutative deformations. 

By the Kostant theorem, \cite{Dix}, $U(\g)$ is a free module over its center.
It follows that at the point $h=0$ the module $(S\g)_{t,0}$ 
is a free module
over the algebra $\I_{t,0}$. One can easily derive from this that
$(S\g)_{t,h}$ is a free module over $\I_{t,h}$.

Now, let the maximal semisimple orbit $M$ be defined by invariant elements from $\I$.
Consider a character defined by $M$, the algebra homomorphism 
$\lambda:\I\to \C$ 
which takes each element from $\I$ to its value on $M$. 
Then, $\C$ may be considered as an $\I$-module, and the function algebra $A$ 
on $M$ is equal to $S\g/\Ker(\lambda)S\g=S\g\ot_\I\C$.
Extend the character $\lambda$ up to a character 
$\lambda_{t,h}:\I_{h,t}\to\C[[h]][t]$ in the trivial way and consider
$\C[[h]][t]$ as a $\I_{h,t}$-module. The tensor product over $\I_{t,h}$  
\be{*}
A_{t,h}=(S\g)_{t,h}\ot\C[[h]][t]
\ee{*}
is a $\C[[h]][t]$-algebra. It is a free $\C[[h]][t]$-module, since
$(S\g)_{t,h}$ is a free one over $\I_{t,h}$.  

It is obvious, that $A_{0,0}=A$, $A_{t,0}$ gives a quantization of
the KKS bracket on $M$,
and  $A_{t,h}$ is a quotient algebra of $(S\g)_{t,h}$.
\end{proof}

In a next paper we shall prove that 
the quantization 
$(S\g)_{t,h}$ can be restricted on all semisimple orbits.

\begin{conjecture}\label{con4.1}
Can be the quantization $(S\g)_{t,h}$ restricted  
on all orbits (not necessarily semisimple)? 
\end{conjecture}
As we have seen, the corresponding Poisson brackets are strongly restricted
on all the orbits.

In next Section we consider the $\Uh$ invariant quantizations
on semisimple orbits in $\g^*$ for all simple Lie algebras $\g$.
It turns out that in general, on a given orbit there are  many nonequivalent
quantizations which are not restrictions from a quantization on $\g^*$. From this
point of view, the quantization on maximal orbits described
by Proposition (\ref{prop3.8}) is a distinguished one.

\section{The one and two parameter quantization\\ 
on semisimple orbits in $\g^*$}

\subsection{Pairs of brackets on semisimple orbits}
Let $\g$ be a simple complex Lie algebra, $\hh$ a fixed Cartan subalgebra.
Let $\Omega\subset \hh^*$ be  the system of roots corresponding to $\hh$.
Select a system of positive roots,  $\Omega^+$, and 
denote by $\Pi\subset\Omega$ the subset of simple roots. 
Fix an element $E_\aa\in \g$  of weight $\alpha$ for each $\alpha\in\Omega^+$ 
and choose $E_{-\aa}$ such that
\be{}\label{Ean}
(E_\aa,E_{-\aa})=1 
\ee{}
for the Killing form $(\cdot,\cdot)$ on $\g$.

Let $\Ga$ be a subset of $\Pi$. Denote by $\hh^*_\Ga$
the subspace in $\hh^*$ generated by $\Ga$. 
Note, that $\hh^*=\hh^*_\Ga\oplus\hh^*_{\Pi\setminus\Ga}$,
and one can identify $\hh^*_{\Pi\setminus\Ga}$ and $\hh^*/\hh^*_\Ga$
via the  projection $\hh^*\to\hh^*/\hh^*_\Ga$.

Let $\Omega_\Ga\subset \hh^*_\Ga$ 
be the subsystem of roots
in $\Omega$ generated by $\Ga$, i.e., 
$\Omega_\Ga=\Omega\cap \hh^*_\Ga$. 
Denote by $\g_\Ga$ the subalgebra of $\g$
generated by the elements $\{E_\aa,E_{-\aa}\}$, $\aa\in\Ga$, and $\hh$.
Such a subalgebra is called the Levi subalgebra.

Let $G$ be a complex connected Lie group with Lie algebra $\g$ and  
$G_\Ga$ a subgroup with Lie algebra $\g_\Ga$. Such a subgroup
is called the Levi subgroup. It is known that
$G_\Ga$ is a connected subgroup. 
Let $M$ be a homogeneous space
of $G$ and $G_\Ga$ be the stabilizer of a point $o\in M$. We
can identify $M$ and the coset space $G/G_\Ga$. 
It is known, that such $M$ is isomorphic to a semisimple orbit in $\g^*$.
This orbit goes through an element $\lambda\in\g^*$ which is just the 
trivial extension to all of  $\g^*$ (identifying
$\g$ and $\g^*$ via the Killing form) of a map 
$\lambda:\hh_{\Pi\setminus\Ga}\to\C$ such that $\lambda(\alpha)\neq 0$
for all $\alpha\in \Pi\setminus\Ga$.
Conversely, it is easy to show that any semisimple orbit in $\g^*$ 
is isomorphic to the quotient of $G$ by a Levi subgroup.

The projection $\pi:G\to M$ induces the map
$\pi_*:\g\to T_o$, where $T_o$ is the tangent space to $M$ at the
point $o$. Since the $\ad$-action of $\g_\Ga$ on $\g$ is semisimple,
there exists an $\ad(\g_\Ga)$-invariant subspace, $\mm=\mm_\Ga$, of $\g$
complementary to $\g_\Ga$, and one can identify $T_o$ and $\mm$ by means of
$\pi_*$. It is easy to see that subspace $\mm$ is uniquely defined
and has a basis formed by the elements
$E_\gamma,E_{-\gamma}$, $\gamma\in\Omega^+\setminus\Omega_\Ga$.

Let $v\in\g^{\ot m}$ be a tensor over $\g$. Using
the right and the left actions of $G$ on itself, one can associate with $v$
right and left invariant tensor fields on $G$ denoted by $v^r$ and $v^l$.

We say that a tensor field, $t$, on $G$ is right $G_\Ga$ invariant, if
$t$ is invariant under the right action of $G_\Ga$.
The $G$ equivariant diffeomorphism between $M$ and $G/G_\Ga$ implies that
any right $G_\Ga$ invariant tensor field $t$ on $G$
induces tensor field $\pi_*(t)$ on $M$. The field $\pi_*(t)$ will be
invariant on $M$ if, in addition, $t$ is left invariant on $G$, 
and any invariant
tensor field on $M$ can be obtained in such a way.
Let $v\in\g^{\ot m}$. For $v^l$ to be right $G_\Ga$
invariant it is necessary and sufficient that $v$ to be $\ad(\g_\Ga)$
invariant. Denote $\pi^r(v)=\pi_*(v^r)$ for any tensor $v$ on $\g$
and $\pi^l(v)=\pi_*(v^l)$ for any $\ad(\g_\Ga)$ invariant tensor 
$v$ on $\g$. 
Note, that tensor $\pi^r(v)$ coincides with the image of $v$ by
the map $\g^{\t m}\to{\rm Vect}(M)^{\t m}$ induced by the action map 
$\g\to{\rm Vect}(M)$.
Any $G$ invariant tensor on $M$ has the form
$\pi^l(v)$. Moreover, $v$ clearly can be uniquely chosen from $\mm^{\ot m}$. 

Denote by $\[v,w\]\in\wedge^{k+l-1}\g$ the Schouten bracket 
of the polyvectors $v\in\wedge^k\g$, $w\in\wedge^l\g$, defined
by the formula 
\be{*}
\[X_1\wedge\cdots\wedge X_k, Y_1\wedge\cdots\wedge Y_l\]=\sum
(-1)^{i+j}[X_i,Y_j]\wedge X_1\wedge\cdots \hat X_i \cdots
\hat Y_j\cdots \wedge Y_l,
\ee{*}
where $[\cdot,\cdot]$ is the bracket in $\g$.
The Schouten bracket is defined in the same way for polyvector fields on
a manifold, but instead of $[\cdot,\cdot]$ one uses the Lie bracket of
vector fields. We will use the same notation for the Schouten
bracket on manifolds.
It is easy to see that $\pi^r(\[v,w\])=\[\pi^r(v),\pi^r(w)\]$, and
the same relation is valid for $\pi^l$.

Denote by $\bO$ the image of $\Omega$ in $\hh^*_{\Pi\setminus\Ga}$
without zero.
It is clear that $\Omega_{\Pi\setminus\Ga}$ can be identified with
a subset of $\bO$ and each element
from $\bO$ is a linear combination of elements from $\Pi\setminus\Ga$
with integer coefficients which all are either positive or negative.
Thus, the subset $\bO^+\subset\bO$ of the elements with positive
coefficients is exactly the image of $\Omega^+$.
We call elements of $\bO$ quasiroots and the images  of $\Pi\setminus\Ga$
simple quasiroots.

\begin{propn}\label{prop5.1}
The space $\mm$ considered as
a $\g_\Ga$ representation space decomposes into the direct sum of
subrepresentations $\mm_\bbe$, $\bbe\in \bO$,
where $\mm_\bbe$ is generated by all the elements $E_\beta$, $\beta\in \Omega$,
such that the projection of $\beta$ is equal to $\bbe$. This decomposition
have the following properties:

a) all $\mm_\bbe$ are irreducible;

b) for $\bbe_1,\bbe_2\in \bO$ such
that $\bbe_1+\bbe_2\in \bO$ one has
$[\mm_{\bbe_1},\mm_{\bbe_2}]=\mm_{\bbe_1+\bbe_2}$;

c) for any pair $\bbe_1,\bbe_2\in \bO$  the representation 
$\mm_{\bbe_1}\ot\mm_{\bbe_2}$ is multiplicity free.
\end{propn}

\begin{proof}
Statements a) and b) are proven in \cite{DGS}. Statement c) follows
from the fact that all the weight subspaces for all $\mm_\bbe$ have 
the dimension one (see N.Bourbaki, Groupes et alg\`ebres de Lie,
Chap. 8.9, Ex. 14).
\end{proof}

Since $\g_\Ga$ contains the Cartan subalgebra $\hh$, each $\g_\Ga$ invariant
tensor over $\mm$ has to be of weight zero. It follows that there are
no invariant vectors in $\mm$. Hence, there are no invariant vector
fields on $M$.

Consider the invariant bivector fields on $M$.
{}From the above, such fields correspond to the $\g_\Ga$ invariant
bivectors from $\wedge^2\mm$. Note, that any $\hh$ invariant bivector
from $\wedge^2\mm$ has to be of the form 
$\sum c(\aa)\Ea\wedge\Ema$.

\begin{propn}\label{pr0}
A bivector  $v\in\wedge^2\mm$ is $\g_\Ga$ invariant if and only if
it has the form
$v=\sum c(\aa)\Ea\wedge\Ema$
where the sum runs over $\aa\in \rc$,
and for two roots $\aa, \beta$ which give the same element
in $\hh^*/\hh^*_\Ga$ 
one has $c(\aa)=c(\beta)$.
\end{propn}

\begin{proof} 
Follows from  Proposition \ref{prop5.1} and condition
(\ref{Ean})
\end{proof}

This proposition shows, that coefficients of an invariant element
$v=\sum c(\aa)\Ea\wedge\Ema$ depend only of the image of $\aa$ in $\bO^+$,
denoted $\baa$, so $v$ can be written in the form
$v=\sum c(\baa)\Ea\wedge\Ema$. 
Let $v\in\wedge^2\mm$ be of the form $v=\sum c(\baa)\Ea\wedge\Ema$,
where the sum runs over $\aa\in \rc$. Denote by $\theta$ the Cartan
automorphism of $\g$. Then, $v$ is $\theta$ anti-invariant, i.e., 
$\theta v=-v$. Hence, any $\g_\Ga$ invariant bivector is $\theta$ 
anti-invariant.
If  $v,w\in\wedge^2\mm$ are $\g_\Ga$ invariant, then 
$\[v,w\]$ is $\theta$ invariant and is of the form 
$\[v,w\]=\sum e(\baa,\bbe)E_{\aa+\beta}\wedge E_{-\aa}\wedge E_{-\beta}$
where roots $\aa,\beta$ are both negative or both positive
and $e(\baa,\bbe)=-e(-\baa,-\bbe)$.
Hence, to calculate  $\[v,w\]$ for such $v$ and $w$ 
it is sufficient to calculate
coefficients $e(\baa,\bbe)$ for positive $\baa$ and $\bbe$.  

Define by $\ff_M$ the invariant three-vector field on $M$ determined
by the invariant element $\ff\in\wedge^3\g$.
A direct computation shows (see \cite{DGS}) that the Schouten bracket of
bivector $v=\sum c(\baa)E_{\aa}\wedge E_{-\aa}$ 
with itself is equal to $K^2\ff_M$ for a complex number $K$, if and only if
the following equations hold
\be{}\label{ff} 
c(\baa+\bbe)(c(\baa)+c(\bbe))=c(\baa)c(\bbe)+K^2
\ee{}
for all the pairs of positive quasiroots $\baa, \bbe$ such that 
$\baa+\bbe$ is a quasiroot.
So, if $c(\baa)$ and $c(\bbe)$ are given and $c(\baa)+c(\bbe)\neq 0$,
\be{}\label{ff1}
c(\baa+\bbe)=\frac{c(\baa)c(\bbe)+K^2}{c(\baa)+c(\bbe)}.
\ee{}

\begin{propn}\label{prop5.3}
Let  $(\baa_1,\cdots,\baa_k)$ be the $k$-tuple of all simple quasiroots.
Given a $k$-tuple of complex numbers $(c_1,...,c_k)$, assign to each $\baa_i$
the number $c_i$.
Then

a) for almost all $k$-tuples of complex numbers (except an algebraic subset
in $\C^k$ of lesser dimension) equations
(\ref{ff1}) uniquely define numbers $c(\baa)$ for all 
positive quasiroots $\baa=\sum \baa_i$ such that the bivector
$f=\sum c(\baa)E_{\alpha}\wedge E_{-\alpha}$ satisfies the condition
$$\[f,f\]=K^2\ff_M;$$

b) when $K=0$, the solution  described in part a)
 defines a Poisson bracket on $M$. Numbers $c(\baa)$
give a solution of (\ref{ff}) if and only if 
there exists a linear form $\ll\in \hh^*_{\Pi\setminus\Ga}$ such that
\be{}\label{ll}
c(\baa)=\frac{1}{\ll(\baa)}
\ee{}
for all quasiroots $\baa$.
\end{propn}

\begin{proof} See  \cite{DGS}.
\end{proof}

\begin{rem}\label{rem5.0}

This proposition shows that invariant brackets $f$ on $M$ defined by part a)
of the proposition 
form a $k$-dimensional variety, $\X_K$,
where $k$ is the number of simple quasiroots. On the other hand,
$k=\dim H^2(M)$, \cite{Bo}. 
If $K$ is regarded as indeterminate, then
$f$ forms a $k+1$ dimensional variety, $\X\subset \C^k\times\C$,
(component $\C$ corresponds to $K$).
Subvariety $\X_0$ corresponds to $K=0$, i.e., consists of
Poisson brackets. It is easy to see that all the Poisson brackets
with $c(\baa)=1/\ll(\baa)\neq 0$ are nondegenerate. Since $\X$ is 
connected, 
it follows that almost all brackets $f$ (except an algebraic subset 
in $\X$ of lesser dimension) are nondegenerate as well.
\end{rem}

\begin{rem}\label{rem5.0a}
Equations (\ref{ff1}) show that when $c(\baa)+c(\bbe)=0$, there appears
a harm for determining $c(\baa+\bbe)$ from given $c(\baa)$ and $c(\bbe)$.
Nevertheless, it is easy to derive from equations (\ref{ff}) that

(*) If $c(\baa)+c(\bbe)=0$ then necessarily $c(\baa)=\pm K$, $c(\bbe)=\mp K$.\\
So it is naturally to consider the quasiroots $\baa$ where $c(\baa)$ are equal to $\pm K$ or not
separately.

Let $c(\baa)$, $\baa\in\bO$, be a solution of equations (\ref{ff}) (we assume 
$c(-\baa)=-c(\baa)$).
It is easy to derive from equations (\ref{ff}) the following properties.

(**)\   If $c(\baa)=\pm K$ and $c(\bbe)\neq \pm K$, then $c(\baa+\bbe)=\pm K$ and 
$c(\baa-\bbe)=\pm K$;

(***) If $c(\baa)=\pm K$ and $c(\bbe)=\pm K$, then $c(\baa+\bbe)=\pm K$.\\
Let $\bO^\prime\subset \bO$ be the subset of quasiroots $\baa$ such that
$c(\baa)\neq\pm K$.
From (**) follows that $\bO^\prime$ is a linear subset, i.e., 
$\bO^\prime=\bO\cap {\rm span}(\bO^\prime)$, where 
${\rm span}(\bO^\prime)$ is the vector subspace of 
$\hh^*/\hh^*_\Gamma$ generated by $\bO^\prime$. 
Let  $(\baa_1,\cdots,\baa_k)$ be a $k$-tuple of elements from $\bO^\prime$ that
form a basis of ${\rm span}(\bO^\prime)$.  
Since by (*)  $c(\baa)+c(\bbe)\neq 0$ for any
$\baa,\bbe\in \bO^\prime$, all $c(\baa)$, $\baa\in \bO^\prime$,
can be found from (\ref{ff1}) using the initial values $c_i=c(\baa_i)$, 
as in Proposition \ref{prop5.3}.

Note that since $c_i\neq \pm K$, there are uniquely defined 
complex numbers $\ll_i\neq 0,1$ such that
$c(\baa_i)=c_i=K\psi(\ll_j)$, where 
\be{*}
\psi(x)=\frac{x+1}{x-1}.
\ee{*}
 Using the formula
\be{*}
\psi(xy)=\frac{\psi(x)\psi(y)+1}{\psi(x)+\psi(y)},
\ee{*}
it is easy to derive that if $\ll:\bO^\prime\to \C^*$ is the multiplicative map
(such that if $\baa,\bbe,\baa+\bbe\in\bO^\prime$ then $\ll(\baa+\bbe)=\ll(\baa)\ll(\bbe)$ ) 
defined by $c(\baa_i)=\ll_i$, then the solution of (\ref{ff1}) is given by the formula
\be{}\label{solff1}
c(\baa)=K\psi(\ll(\baa)), \qqquad \baa\in \bO^\prime.
\ee{}
For correctness of this formula, one needs that the map $\ll$ to be regular, i.e.,
that $\ll$ 
to satisfy the condition:
if $\baa,\bbe,\baa+\bbe\in\bO^\prime$
then $\ll(\baa)\ll(\bbe)=1$ only when $\baa=-\bbe$.

From property (**) follows that the numbers $c(\baa)$ define a function on
the set $\pi(\bO)$, where $\pi$ is the natural map 
$\hh^*/\hh^*_\Gamma \to (\hh^*/\hh^*_\Gamma)/ {\rm span}(\bO^\prime)$.
This function has values $\pm K$. Let $X\subset\pi(\bO)$ be the subset where
this function has value $K$.
From property (***) follows that $X$ is a semilinear subset. It means that
if $x_1,x_2\in X$ and $x_1+x_2\in \pi(\bO)$ then $x_1+x_2\in X$, and
$X\cap(-X)=\emptyset$, $X\cup(-X)=\pi(\bO)$.

The arguments above lead to the following description of the variety $\Z_K$ of all
solutions of (\ref{ff}) (or, what is the same, the variety of invariant
brackets $f$ on $M$ such that $\[f,f\]=K^2\ff_M$).

\begin{propn}\label{prop5.3a}
Variety $\Z_K$ splits into stratas. Each strata is defined by choosing a linear
subset $\bO^\prime$ of $\bO$ and a semilinear subset $X$ of $\pi(\bO)$. 
Points of this strata are parameterized by the
multiplicative regular maps $\ll:\bO^\prime\to \C^*$.

Let the data $(\bO^\prime, X, \ll)$ corresponds to a point of $\Z_K$.
Then the coefficients $c(\baa)$ of $f$ are determined in the following way.
If $\baa\in \bO^\prime$ then $c(\baa)$ is found by (\ref{solff1}).
If $\pi(\baa)\in X$ then $c(\baa)=K$.
If $\pi(\baa)\in -X$ then $c(\baa)=-K$.
\end{propn}
Of course, in case $K=0$ the choose of $X$ does not matter: a strata of $\Z_0$ is
determined only by choosing   $\bO^\prime$.

Note also that the description of $\Z_K$ given in the proposition does not
depend on choosing a basis in $\bO$. The variety $\X_K$ from the previous remark
forms an open everywhere dense subset of $\Z_K$ and does depend on choosing
a basis. According to Remark \ref{remPL} this proposition describes all the
$(G,\tilde{r})$-Poisson structures on semisimple orbits.

\end{rem}

%**********************************************************************

Now we fix a Poisson bracket
$s=\sum (1/\ll(\baa))E_{\aa}\wedge E_{-\aa}$, 
where $\ll$ is a fixed linear form,
and
describe the invariant brackets 
$f=\sum c(\baa)E_{\aa}\wedge E_{-\aa}$ which satisfy 
the conditions
\be{}\label{pbf}
&\[f,f\]=K^2\ff_M  \qqquad {\rm for} \quad K\neq 0, \\ 
&\[f,s\]=0. \notag
\ee{}

Direct computation shows that the condition $\[f,s\]=0$ is equivalent 
to the system of equations for the coefficients $c(\baa)$ of $f$
\be{} \label{comp}
c(\baa)\ll(\baa)^2+c(\bbe)\ll(\bbe)^2=c(\baa+\bbe)\ll(\baa+\bbe)^2
\ee{}
for all the pairs of positive quasiroots $\baa, \bbe$ such that 
$\baa+\bbe$ is a quasiroot.

\begin{defn}
Let $M$ be an orbit in $\g^*$ (not necessarily semisimple).
We call $M$ a {\it good} orbit, if there exists an invariant bracket,
$f$,  on $M$ satisfying 
the conditions (\ref{pbf}) for $s$ the Kirillov-Kostant-Souriau (KKS)
Poisson bracket on $M$.
\end{defn}

So, a semisimple orbit $M$ is a good orbit if and only if 
equations (\ref{ff}) and (\ref{comp})
are compatible, i.e., have a common solution.

\begin{propn}\label{prop5.4} 
The good semisimple orbits are the following:

a) For $\g$ of type $A_n$ all semisimple orbits are good.

b) For  all other $\g$, the orbit $M$ is good if and only if
the set $\Pi\setminus\Ga$ consists of one or two roots which appear   
in representation of the maximal root with coefficient 1.

c) The brackets $f$ on good orbits form a one-dimensional variety:
all such brackets have the form 
\be{*}
\pm f_0+ts,
\ee{*}
where $t\in \C$ and $f_0$ is a fixed bracket satisfying (\ref{pbf}).
\end{propn}

\begin{proof} See \cite{DGS}.
\end{proof}

\begin{rem}\label{rem5.1}
From Proposition (\ref{prop3.5}) follows that for $\g=sl(n)$ all
orbits (not only semisimple) are good ones.
In addition, if an orbit, $M$, is such that $\ff_M=0$, then $M$ is good:
one can take $f=0$.
In \cite{GP} there is a classification of orbits for all simple $\g$, for which
$\ff_M=0$.
\end{rem}

\begin{conjecture}\label{con5.1}
Let $\g$ be a simple Lie algebra. Are all orbits in $\g^*$ good?
If not, what is a classification of good orbits?
\end{conjecture} 

\subsection{Cohomologies defined by invariant brackets}

In the next subsection we prove the existence of a 
$U_h(\g)$ invariant quantiztion 
of the Poisson brackets described above using the methods 
of \cite{DS1}. This requires us to consider the $3$-cohomology of
the complex 
$(\Lambda^\bullet (\g/\g_\Ga))^{\g_\Ga}=(\Lambda^\bullet\mm)^{\g_\Ga}$
of $\g_\Ga$ invariants
with differential given by the Schouten
bracket with the bivector 
$f\in(\Lambda^2\mm)^{\g_\Ga}$ from Proposition 
\ref{prop5.3} a), 
$$\delta_f:u\mapsto \[f,u\]\,\qqquad\mbox{for} 
\quad u\in(\Lambda^\bullet\mm)^{\g_\Ga}.$$ 
The condition  $\delta_f^2=0$ follows from the Jacobi identity for
the Schouten bracket together with the fact that $\[f,f\]=K^2\ff_M$.
Denote these cohomologies by $H^k(M,\delta_f)$, whereas the usual
de Rham cohomologies are denoted by $H^k(M)$.
 
Recall (see Remark \ref{rem5.0}) that the brackets $f$ satisfying $\[f,f\]=K^2\ff$ 
form a connected variety $\X$ which contains a submanifold $\X_0$
of Poisson brackets.

\begin{propn}\label{prop5.5} For almost all $f\in \X$ 
(except an algebraic subset
of lesser dimension) one has
\be{*}
H^k(M,\delta_f)=H^k(M)
\ee{*}
for all $k$.
In particular, $H^k(M,\delta_f)=0$ for odd $k$.
\end{propn}

\begin{proof}
First, let $v$ be a Poisson bracket, i.e., $v\in\X_0$.
Then the complex of polyvector fields on $M$, $\Theta^\bullet$, with
the differential $\delta_v$ is well defined.
Denote by $\Omega^\bullet$ the de Rham complex on $M$.
Since none of the coefficients $c(\baa)$ of $v$ are
zero, $v$ is a nondegenerate bivector field, and therefore
it defines an $\A$-linear isomorphism 
$\tilde{v}:\Omega^1\to\Theta^1$, $\omega\mapsto v(\omega,\cdot)$,
which can be extended up to the isomorphism $\tilde{v}:\Omega^k\to\Theta^k$
of $k$-forms onto $k$-vector fields for all $k$.
Using Jacobi identity for $v$ and invariance of $v$, one can
show that $\tilde{v}$ gives a $G$ invariant isomorphism of
these complexes, so their cohomologies are the same.

Since $\g$ is simple, the subcomplex of 
$\g$ invariants, $(\Omega^\bullet)^\g$, splits off
as a subcomplex of $\Omega^\bullet$.
In addition, $\g$ acts trivially on cohomologies,
since for any $g\in G$ the map $M \to M$,
$x\mapsto gx$, is homotopic to the identity map, ($G$ is a connected
Lie group corresponding to $\g$).
It follows that cohomologies of complexes
$(\Omega^\bullet)^\g$ and $\Omega^\bullet$ coincide.

But $\tilde{v}$ gives an isomorphism of complexes
$(\Omega^\bullet)^\g$ and 
$(\Theta^\bullet)^\g=((\Lambda^\bullet\mm)^{\g_\Ga},\delta_v)$.
So, cohomologies of the latter complex coincide with de Rham
cohomologies, which
proves the proposition for $v$ being Poisson brackets.

Now, consider the family of complexes
$((\Lambda^\bullet\mm)^{\g_\Ga},\delta_v)$, $v\in\X$. It is clear that
$\delta_v$ depends algebraicly on $v$. It follows from the uppersemicontinuity  
of $\dim H^k(M,\delta_v)$ and the fact that $H^k(M)=0$ for odd $k$,
\cite{Bo}, that $H^k(M,\delta_v)=0$ for odd $k$ and almost all $v\in\X$.
Using the uppersemicontinuity again and the fact that
the number $\sum_k(-1)^k\dim H^k(M,\delta_v)$ is the same for all $v\in\X$,
we conclude that $\dim H^k(M,\delta_v)=\dim H^k(M)$ for even $k$ and
almost all $v$.
\end{proof}  

\begin{rem}
Call $f\in \X$ admissible, if it satisfies Proposition \ref{prop5.5}.
{}From  the proof of the proposition follows that the subset $\D$ 
such that $\X\setminus\D$ consists of admissible brackets
does not intersect with the subset $\X_0$
consisting of Poisson brackets. 
It follows from this fact that for  each good orbit there are admissible $f$ 
compatible with the KKS bracket.
Indeed,
let $M$ be a good orbit and $f_0+ts$ the family from Proposition \ref{prop5.4} c)
satisfying (\ref{pbf}) for a fixed $K$. Then for almost all numbers $t$ this
bracket is admissible. In fact, this family is contained in the two parameter
family $uf_0+ts$. By $u=0$, $t\neq 0$ we obtain admissible brackets.
So, there exist $u_0\neq 0$ and $t_0$ such that the bracket
$u_0f_0+t_0s$ is admissible. It follows that the bracket
$f_0+(t_0/u_0)w$ is admissible, too. So, in the family $f_0+ts$
there is an admissible bracket, and we conclude that almost all brackets
in this family (except a finitely many) are admissible.
\end{rem}

For the proof of existence of two parameter quantization
for the cases $D_n$ and $E_6$ in the next subsection,
we will use the following result on invariant three-vector fields.

Denote by $\theta$ the Cartan automorphism of $\g$.

\begin{lemma}\label{lem5.1}
For either $D_n$ or $E_6$ and one of the
subsets, $\Ga$, of simple roots such that $G_\Ga$ defines a good orbit, 
any  $\g_\Ga$ and $\theta$
invariant element $v$ in $\Lambda^3\mm$ is a multiple of $\ff_M$, that is,
 $$\left(\Lambda^3(\mm\right)^{\g_\Ga}\cong\langle\ff_M\rangle.$$ 
\end{lemma}     

\begin{proof}
In this case the system of positive quasiroots consists of $\baa$, $\bbe$, and $\baa+\bbe$,
where $\baa$, $\bbe$ are the simple quasiroots. From Proposition \ref{prop5.1}
follows that invariant elements 
in $\mm_{\baa}\t \mm_{\bbe}\t \mm_{-\baa-\bbe}$ and
$\mm_{-\baa}\t \mm_{-\bbe}\t \mm_{\baa+\bbe}$ form 
subspaces of dimension one, $I_1$ and $I_2$.
Moreover, all the ivariant elements of $\Lambda^3\mm$ are lying in
$I_1+I_2$. Since $\theta$ takes $I_1$ onto $I_2$, there is only one-dimansional
$\theta$ invariant subspace in $I_1+I_2$, which is necessarily
generated by $\ff_M$.
\end{proof}

\subsection{$\Uh$ invariant quantizations in one and two parameters}

In this subsection we prove the existence of one and two parameter
$\Uh$ invariant quantization
of the function algebras $\A$ on  semisimple orbits, $M$, in $\g^*$. 
By Proposition \ref{prop2.2}, 
the one parameter quantization has the Poisson bracket of the form
\be{}\label{u1}
f(a,b)-\{a,b\}_r, \qquad \[f,f\]=-\ff_M.
\ee{}
We show that the one parameter quantization exists
for all semisimple orbits and all $f$ constructed in
Proposition \ref{prop5.3} a) and satisfying
Proposition \ref{prop5.5}. 

For two parameter quantization, there are two compatible Poisson brackets:
the KKS bracket $s$ and
the bracket of the form (\ref{u1}) with the additional condition
\be{}\label{u2}
\[f,s\]=0.
\ee{}
We show that the two parameter quantization exists for good orbits 
in cases $D_n$ and $E_6$ and for almost all $f$ satisfying (\ref{u1})
and (\ref{u2}).

Note that in subsection \ref{ss3.5} we have proven that in case $A_n$
the two parameter quantization exists for maximal semisimple orbits.
In a next paper we shall prove the same for all semisimple orbits.

We remind the reader of the method in  \cite{DS1}.
The first step is to construct a $U(\g)$ invariant
quantization in the category ${\cA}(U(\g)[[h]],\De,\Phi_h)$.
Then we use the equivalence
given by the pair $(\Id,F_h)$  between the monoidal categories 
${\cA}(U(\g)[[h]],\De,\Phi_h)$
and ${\cA}(U(\g)[[h]],\De_h,{\bf 1})$ 
to define a $\Uh$ invariant quantization,
either $\mu_h F^{-1}_h$ in the one parameter case or 
$\mu_{t,h}F^{-1}_h$ in the two parameter case (see Subsections 2.2 and 2.3).
In the following we often write $\Phi$ for $\Phi_h$.

\begin{propn}\label{prop5.6}
Let $\g$ be a simple Lie algebra, $M$ a semisimple orbit in $\g^*$.
Then, for almost all (in sense of Proposition \ref{prop5.5})
$\g$ invariant brackets $f$
satisfying $\[f,f\]=-\ff_M$, there exists a 
multiplication  $\mu_h$ on $\A$ 
$$\mu_h(a,b)=ab +(h/2) f(a,b) +\sum_{n\geq 2}h^n \mu_n(a,b)$$
which  is $\Ug$ invariant (equation \ref{finv})) and
$\Phi$ associative (equation (\ref{fass})).
\end{propn}
\begin{proof} To begin,  consider the multiplication $\mu^{(1)}(a,b)=ab +(h/2) f(a,b)$.
The corresponding  obstruction cocycle  is given by 
$$obs_2=\frac1{h^2}(\mu^{(1)}(\mu^{(1)}\otimes id)-
\mu^{(1)}(id\otimes\mu^{(1)})\Phi)$$
considered modulo terms of order $h$. No $\frac1h$ terms appear because
$f$ is a biderivation and, therefore, a Hochschild cocycle. 
The fact that the presence of
$\Phi$ does not interfere with the cocyle condition and 
that this equation defines a Hochschild $3$-cocycle 
was proven in \cite{DS1}. It is well known that if we restrict to
the subcomplex of cochains given by differential operators, the
differential Hochschild cohomology of $\A$ in dimension $p$ is 
the space of $p$-polyvector fields on $M$. 
Since $\g$ is reductive, the subspace of $\g$ invariants splits off
as a subcomplex and has cohomology given by $(\Lambda^p\mm)^{\g_\Ga}$.
The complete antisymmetrization of a $p$-tensor projects
the  space of invariant differential $p$-cocycles onto the subspace
$(\Lambda^p\mm)^{\g_\Ga}$ representing the cohomology. The equation
$\[f,f\]+\ff_M=0$ implies that obstruction cocycle is a coboundary,
and we can find a $2$-cochain $\mu_2$, so that 
$\mu^{(2)}=\mu^{(1)}+ h^2\mu_2$ satisfies 
$$\mu^{(2)}(\mu^{(2)}\otimes id)-
\mu^{(2)}(id\otimes\mu^{(2)})\Phi=0\mbox{  mod }h^2.$$
Assume we have defined the deformation $\mu^{(n)}$ to order $h^n$ 
such that $\Phi$ associativity holds modulo $h^n$, then we define the
$(n+1)^{\mbox{st}}$ obstruction cocycle by
$$obs_{n+1}=\frac1{h^{n+1}}(\mu^{(n)}(\mu^{(n)}\otimes id)-
\mu^{(n)}(id\otimes\mu^{(n)})\Phi)\mbox{ mod }h.$$

In \cite{DS1} (Proposition 4) we showed that the usual 
proof that the obstruction cochain satisfies the cocycle condition
carries through to the $\Phi$ associative case. 
The coboundary of $obs_{n+1}$ appears as the 
$h^{n+1}$ coefficient of the signed sum of the compositions of 
$\mu^{(n+1)}$ with $obs_{n+1}$. 
The fact that  $\Phi=1$ mod $h^2$ together with  the pentagon
identity implies that the sum vanishes
identically, and thus all coefficients vanish, including the coboundary
in question. Let $obs_{n+1}'\in (\Lambda^3\mm)^{\g_\Ga}$ be 
the projection of $obs_{n+1}$ on the totally skew symmetric part, which represents
the cohomology class of the obstruction cocycle. The coefficient of
$h^{n+2}$ in the same signed sum, when projected on the skew symmetric
part, is $\[f,obs_{n+1}'\]$ which is the coboundary of $obs'_{n+1}$  
in the complex $(\Lambda^\bullet\mm)^{\g_\Ga},\delta_f= \[f,.\])$. Thus
$obs'_{n+1}$  is a $\delta_f$ cocycle.  By Proposition \ref{prop5.5},
this complex has zero cohomology. Now we  modify 
$\mu^{(n+1)}$ by adding a term $h^n\mu_n$ with $\mu_n\in
(\Lambda^2\mm)^{\g_\Ga}$ and consider the $(n+1)^{\mbox{st}}$ obstruction 
cocycle for $\mu^{\prime(n+1)}=\mu^{(n+1)}+h^n\mu_n$. Since the term we added at
degree $h^n$ is a Hochschild cocyle, we do not introduce a $h^n$ term
in the calculation of $\mu^{(n)}(\mu^{(n)}\otimes id)-
\mu^{(n)}(id\otimes\mu^{(n)})\Phi $ and the totally skew symmetric
projection $h^{n+1}$ term has been  modified
by $\[f,\mu_n\]$. By choosing $\mu_n$ appropriately,
we can make the $(n+1)^{\mbox {st}}$ obstruction cocycle represent the
 zero cohomology class, and we are able to continue
the recursive construction of the desired deformation.
\end{proof}
  
Now we prove the existence of a two parameter deformation 
for good orbits in the cases $D_n$ and $E_6$.
\begin{propn}\label{prop5.7}
Given a pair of $\g$ invariant brackets, $f,v$, on
a good orbit in $D_n$ or $E_6$ 
satisfying $\[f,f\]=-\ff_M$, $\[f,v\]=\[v,v\]=0$, there exists a 
multiplication  $\mu_{h,t}$ on $\A$ 
$$\mu_{t,h}(a,b)=ab+(h/2)f(a,b)+(t/2)v(a,b)+
\sum_{k,l\geq 1} h^kt^l\mu_{k,l}(a,b)$$
which  is $\Ug$ invariant and
$\Phi$ associative.
\end{propn}
\begin{proof} The existence of a multiplication which is $\Phi$ associative
up to and including $h^2$ terms is nearly identical to the previous 
proof. Both $f$ and $v$ are anti-invariant under the Cartan involution
$\theta$. 
We shall look for a multiplication  $\mu_{t,h}$ such that
$\mu_{k,l}$ is $\theta$ anti-invariant and skew-symmetric for odd $k+l$
and $\theta$ invariant and symmetric for even $k+l$.

 So, suppose we have a multiplication 
defined to order $n$,
$$\mu_{t,h}(a,b)=ab+h\mu_1(a,b)+t\mu^\prime_1(a,b)+
\sum_{k+l\leq n} h^kt^l\mu_{k,l}(a,b),$$
with mentioned above invariance properties
and $\Phi$ associative to order $h^n$.

Further we shall suppose that $\Phi$ has the properties:
It  is invariant under the Cartan involution $\theta$ and $\Phi^{-1}=\Phi_{321}$.
Such $\Phi$ always can be choosen, \cite{DS2}. %coh constr of quan groups
Using these properties for $\Phi$,
direct computation shows that the obstruction cochain, 
$$obs_{n+1}=\sum_{k=0,\ldots, n+1} h^k t^{n+1-k}\beta_k,$$
has the following invariance properties:
For odd $n$, $obs_{n+1}$ is $\theta$ invariant and 
$obs_{n+1}(a,b,c)=-obs_{n+1}(c,b,a)$, and for even $n$,  
and $obs_{n+1}$ is $\theta$ anti-invariant and 
$obs_{n+1}(a,b,c)=obs_{n+1}(c,b,a)$.

Hence, the projection of $obs_{n+1}$ on $(\Lambda^3\mm)^{\g_\Ga}$ is
equal to zero for even $n$. It follows that all the $\beta_k$ are 
Hochschild coboundaries, and the standard argument implies
that the multiplication can be extended up to order $n+1$ with
the required properties. 

For odd $n$, Lemma \ref{lem5.1} shows that
the projection on $(\Lambda^3\mm)^{\g_\Ga}$ has the form
$$obs_{n+1}=\left(\sum_{k=0,\ldots, n+1} a_kh^k t^{n+1-k}\right)\ff_M.$$
The KKS bracket is given by the
two-vector
$$v=\sum _{\alpha\in\Omega^+\setminus\Omega_\Ga}\frac1{\lambda(\bar\alpha)}
E_\alpha\wedge E_{-\alpha}.$$
Setting 
$$w=\sum _{\alpha\in\Omega^+\setminus\Omega_\Ga}\lambda(\bar\alpha)
E_\alpha\wedge E_{-\alpha},$$
gives
$$\[ v,w\]=-3\ff_M.$$
Defining 
$$\mu^{\prime(n)}=\mu^{(n)}+ \frac{a_0}3 t^n w,$$ 
the new obstruction cohomology class is
$$obs'_{n+1}=(\sum_{k=1,\ldots, n+1} a_kh^k t^{n+1-k})\ff_M.$$
Finally we define 
$$\mu^{\prime\prime(n)}=\mu^{\prime(n)}+\sum_{k=1,\ldots, n+1} a_kh^{k-1} t^{n+1-k})f$$ 
and get an obstruction cocycle which is zero in cohomology.
Now the standard argument implies that the deformation can be extended
to give a $\Phi$ associative invariant multiplication with the required
properties
of order $n+1$.

So, we are able to continue the recursive construction of the
desired multiplication. 
\end{proof} 

Using the $\Phi_h$ associative multiplications $\mu_h$ and $\mu_{t,h}$ 
from Propositions \ref{prop5.6} and \ref{prop5.7}
and the equivalence between the monoidal categories 
${\cA}(U(\g)[[h]],\De,\Phi_h)$
and ${\cA}(U(\g)[[h]],\wt{\De},{\bf 1})$ 
given by the pair $(\Id,F_h)$ (see Section 2), 
one can define $\Uh$ invariant multiplications,
either $\mu_h F^{-1}_h$ in the one parameter case or 
$\mu_{t,h}F^{-1}_h$ in the two parameter case.

\begin{rem}
After \cite{Ko}, the philosophy is that there are no obstructions for
quantizations of Poisson brackets on manifolds. In this connection, the following question arises:
\begin{conjecture}
Let $M$ be a $G$-manifold on which there exists an invariant connection.
Given a $G$ invariant Poisson bracket, $v$, on $M$, 
does there exist a $G$ invariant quantization of $v$?
\end{conjecture}
In case $M$ is a homogeneous manifold the bracket $v$ has a constant rank,
and such a quantization can be obtained by Fedosov's method, \cite{Fed}, \cite{Do1}.

Another question which relates to the topic of this paper is the following.
\begin{conjecture}
Let $M$ be a $G$-manifold on which there exists an invariant connection,
$\Ug$ the corresponding to $G$ universal enveloping algebra,
and $\Phi_h\in(\Ug)^{\ot 3}[[h]]$ an invariant element of
the form (\ref{fPhi}) obeying
the pentagon identity (\ref{pent}). 
Let $f$ be an invariant bracket on $M$ 
satisfying $\[f,f\]=-\ff_M$.
Does there exist a $\Ug$ (or $G$) invariant
and $\Phi_h$ associative quantization of $f$ (as in Proposition \ref{prop5.6})?
\end{conjecture}
Note that if the answer to this Question is positive, then the answer to
Question \ref{con2.1} is also positive:  we take for $M$ the group $G$ itself
and consider it as a $G$-manifold by left multiplication.

\end{rem}

\end{document}